\documentclass[ps,dvips,noinfoline]{imsart}
\RequirePackage{hyperref}

\doi{10.1214/07-PS122}
\pubyear{2007}
\volume{4}
\firstpage{268}
\lastpage{302}

\startlocaldefs
\usepackage{amssymb,epsfig} 
\newtheorem{lemma}{Lemma}[section]
\newtheorem{definition}{Definition}[section]
\newtheorem{theorem}{Theorem}[section]
\newtheorem{coro}{Corollary}[section]

\newtheorem{remark}{Remark}[section]

\newcommand{\bsq}{\vrule height .9ex width .8ex depth -.1ex}
\newcommand{\RR}{{\mathbb R}}

\newcommand{\sD}{{\cal D}}
\newcommand{\sC}{{\cal C}}

\newcommand{\sF}{{\cal F}}

\newcommand{\sS}{{\cal S}}

\newcommand{\hsp}{\hspace*{\parindent}}
\newcommand{\eeq}{\end{equation}}

\newcommand{\ra}{\rightarrow}

\newcommand{\beql}[1]{\begin{equation}\label{#1}}
\newcommand{\eqn}[1]{(\ref{#1})}
\newcommand{\beq}{\begin{displaymath}}
\newcommand{\eeqno}{\end{displaymath}}

\newcommand{\qandq}{\quad\mbox{and}\quad}
\newcommand{\qforq}{\quad\mbox{for}\quad}

\newcommand{\qasq}{\quad\mbox{as}\quad}
\newcommand{\qinq}{\quad\mbox{in}\quad}

\newcommand{\qforallq}{\quad\mbox{for all}\quad}

\endlocaldefs

\begin{document}

\begin{frontmatter}

\title{Proofs of the martingale FCLT\protect\thanksref{T1}}
\thankstext{T1}{This is an original survey paper}
\runtitle{Martingale FCLT}

\begin{aug}
 \author{\fnms{Ward} \snm{Whitt}\ead[label=e1]{ww2040@columbia.edu}\ead[label=e2,url]{www.columbia.edu/$\sim$ww2040}}
\address{Department of Industrial Engineering and Operations Research\\ Columbia University, New York, NY 10027-6699\\ \printead{e1}\\ \printead{e2}}



\runauthor{W. Whitt}
\end{aug}

\begin{abstract}
This is an expository review paper elaborating on the
proof of the martingale functional central limit theorem (FCLT).
This paper also reviews tightness and stochastic boundedness,
highlighting one-dimensional criteria for tightness used in the proof of the martingale FCLT.
This paper supplements the expository review paper Pang,
Talreja and Whitt (2007) illustrating the ``martingale method''
for proving many-server heavy-traffic stochastic-process limits
for queueing models, supporting diffusion-process approximations.
\end{abstract}



 \received{\smonth{12} \syear{2007}}

\tableofcontents

\end{frontmatter}

\section{Introduction}\label{secIntro}

In this paper we
elaborate upon the proof of the martingale FCLT in \S 7.1 of Ethier
and Kurtz \cite{EK86}.  We also discuss alternative arguments, as in
Jacod and Shiryaev \cite{JS87}.  The proof in \cite{EK86} is correct,
but it is concise, tending to require knowledge of previous parts of the book.
We aim to make the results more accessible by sacrificing some generality and focusing on the ``common case.''

In addition to reviewing proofs of
 the martingale FCLT,
we review tightness criteria, which play an important role in the
proofs of the martingale FCLT.
  An important role is played by simple ``one-dimensional'' criteria for tightness
of stochastic processes in $D$, going beyond the classical criteria in
Billingsley \cite{B68, B99}, reviewed in Whitt \cite{W02}.  The alternative
one-dimensional criteria come from Billingsley \cite{B74, B99}, Kurtz
\cite{K75, K81}, Aldous \cite{A78, A89}, Rebolledo \cite{R80} and Jacod et
al. \cite{JMM83}.\looseness=-1

The martingale FCLT has many applications, so that this paper may aid
in many contexts.  However, we were motivated by applications to queueing models.
Specifically, this paper is intended to
 supplement Pang, Talreja and Whitt \cite{PTW07},
which is an expository review paper illustrating how to
do martingale proofs of many-server heavy-traffic limit theorems
for Markovian queueing models, as in Krichagina and Puhalskii \cite{KP97} and
Puhalskii and Reiman \cite{PR00}.  Pang et al. \cite{PTW07} review martingale basics,
indicate how martingales arise in the queueing models and show
how they can be applied to establish the stochastic-process limits
for the queueing models.  They apply the martingale method to
the elementary $M/M/\infty$ queueing model and a few variations.

The rest of this paper is organized as follows: We start in \S
\ref{secMartFCLT} by stating a version of
 the martingale FCLT from p 339 of Ethier and Kurtz \cite{EK86}.
Next in \S \ref{secTightAll} we review tightness, focusing
especially on criteria for tightness. Then we turn to the proof of
the martingale FCLT. We give proofs of tightness in \S
\ref{secTightProof}; we give proofs of the characterization in \S
\ref{secCharProof}.

\section{The Martingale FCLT}\label{secMartFCLT}

We now state a version of the martingale FCLT for a sequence of
local martingales $\{M_n: n \ge 1\}$ in $D^k$, based on Theorem
7.1 on p. 339 of Ethier and Kurtz \cite{EK86}, hereafter referred to as
EK. We shall also make frequent reference to Jacod and Shiryayev
\cite{JS87}, hereafter referred to as JS. See Section VIII.3 of JS for
related results; see other sections of JS for generalizations.

We will state a special case of Theorem 7.1 of EK
in which the limit process is multi-dimensional Brownian
motion.  However, the framework always produces limits with
continuous sample paths and independent Gaussian increments. Most
applications involve convergence to Brownian motion.  Other
situations are covered by JS, from which we see that proving
convergence to discontinuous processes evidently is more
complicated.

We assume familiarity with martingales; basic notions are reviewed in \S 3
of Pang et al. \cite{PTW07}.
The key part of each condition in the martingale FCLT below is the convergence of the
quadratic covariation processes.  Condition (i) involves the
optional quadratic-covariation (square-bracket) processes
$\left[M_{n,i}, M_{n,j}\right]$, while condition (ii) involves the
predictable quadratic-covariation (angle-bracket) processes
$\langle M_{n,i}, M_{n,j} \rangle$.  Recall
that the square-bracket process is more general, being well
defined for any local martingale (and thus any martingale),
whereas the associated angle-bracket process is well defined only
for any locally square-integrable martingale (and thus any
square-integrable martingale); see \S 3.2 of \cite{PTW07}.

Thus the key conditions below are the assumed convergence of the
quadratic-variation processes in conditions \eqn{ek2} and
\eqn{ek6}. The other conditions \eqn{ek1}, \eqn{ek4} and \eqn{ek6}
are technical regularity conditions. There is some variation in
the literature concerning the extra technical regularity
conditions; e.g., see Rebolledo \cite{R80} and JS \cite{JS87}.

Let $\Rightarrow$ denote convergence in distribution and let
$D \equiv D([0,\infty), \RR)$ be the usual space of right-continuous
real-valued functions on the semi-infinite interval $[0, \infty)$
with limits from the left, endowed with the Skorohod \cite{S56} $J_1$ topology;
see Billingsley \cite{B68, B99}, EK, JS and Whitt \cite{W02} for
background.
For a
function $x$ in $D$, let $J(x,T)$ be the absolute value of the maximum
jump in $x$ over the interval $[0,T]$, i.e.,
 \beql{mod3}
J(x,T) \equiv \sup{\{|x(t) - x(t-)|: 0 < t \le T \}} ~.
\eeq

\begin{theorem}{\em $($multidimensional martingale FCLT$)$}\label{thMart}
For $n \ge 1$, let $M_n \equiv (M_{n,1}, \ldots , M_{n,k})$ be a
local martingale in $D^k$ with respect to a filtration
$\textbf{F}_n \equiv \{\sF_{n,t}: t \ge 0\}$ satisfying $M_n (0) =
(0, \ldots , 0)$.  Let $C \equiv (c_{i,j})$ be a $k \times k$
covariance matrix, i.e., a nonnegative-definite symmetric matrix
of real numbers.

\vspace{0.1in}
  \textbf{Assume that one of the following two conditions holds:}

\vspace{0.1in}
\textbf{$(i)$}  The expected value of the maximum
jump in $M_n$ is asymptotically negligible; i.e., for each
$T > 0$,
\beql{ek1}
 \lim_{n \ra \infty}{\left\{ E\left[
J(M_n,T)\right]\right\}} = 0
\eeq
and, for each pair $(i,j)$ with
$1 \le i \le k$ and $1 \le j \le k$, and each $t > 0$,
\beql{ek2}
\left[M_{n,i}, M_{n,j} \right] (t) \Rightarrow c_{i,j} t \qinq \RR
\qasq n \ra \infty ~.
\eeq

\vspace{0.1in} \textbf{$(ii)$}  The local martingale
$M_n$ is locally square-integrable, so that the predictable
quadratic-covariation processes $\langle M_{n,i}, M_{n,j} \rangle$
can be defined. The expected value of the maximum jump in $\langle
M_{n,i}, M_{n,j} \rangle$ and the maximum squared jump of $M_n$
are asymptotically negligible; i.e., for each $T > 0$ and $(i,j)$
with $1 \le i \le k$ and $1 \le j \le k$, \beql{ek4} \lim_{n \ra
\infty}{\left\{ E\left[ J \left( \langle M_{n,i}, M_{n,j} \rangle,
T \right) \right]\right\}}  = 0 ~, \eeq \beql{ek5} \lim_{n \ra
\infty}{\left\{ E\left[ J \left( M_n, T \right)^2 \right]\right\}}
= 0 ~, \eeq and \beql{ek6} \langle M_{n,i}, M_{n,j} \rangle (t)
\Rightarrow c_{i,j} t \qinq \RR \qasq n \ra \infty \eeq for each
$t > 0$ and for each $(i,j)$.

\vspace{0.1in}
  \textbf{Conclusion:}

\vspace{0.1in} If indeed one of the the conditions $(i)$ or $(ii)$
above holds, then \beql{ek7} M_n \Rightarrow  M \qinq D^k \qasq n
\ra \infty ~, \eeq where $M$ is a \textbf{$k$-dimensional
$(0,C)$-Brownian motion}, having mean vector and covariance matrix
\beql{ek8} E[M(t)] = (0, \ldots , 0) \qandq E[M(t)M(t)^{tr}] = C
t, \quad t \ge 0 ~, \eeq where, for a matrix $A$, $A^{tr}$ is the
transpose.
\end{theorem}

Of course, a common simple case arises when $C$ is a diagonal
matrix; then the $k$ component marginal one-dimensional Brownian
motions are independent. When $C = I$, the identity matrix, $M$ is
a standard $k$-dimensional Brownian motion, with independent
one-dimensional standard Brownian motions as marginals.

At a high level, Theorem \ref{thMart} says that, under regularity
conditions, convergence of martingales in $D$ is implied by
convergence of the associated quadratic covariation processes.  At
first glance, the result seems even stronger, because we need
convergence of only the one-dimensional quadratic covariation
processes for a single time argument. However, that is misleading,
because the stronger weak convergence of these quadratic
covariation processes in $D^{k^2}$ is actually equivalent to the
weaker required convergence in $\RR$ for each $t, i, j$ in
conditions \eqn{ek2} and \eqn{ek6}, as we will show below.

To state the result, let $[[M_n]]$ be the matrix-valued random
element of $D^{k^2}$ with $(i,j)^{\rm th}$ component $[M_{n,i},
M_{n,j}]$; and let $\langle \langle M_n \rangle \rangle$ be the
matrix-valued random element of $D^{k^2}$ with $(i,j)^{\rm th}$
component $\langle M_{n,i}, M_{n,j}\rangle$. Let $e$ be the
identity map, so that $Ce$ is the matrix-valued deterministic
function in $D^{k^2}$ with elements $\{c_{i,j} t: t \ge 0\}$.

\begin{lemma}{$($modes of convergence for quadratic covariation processes$)$}\label{lemModQCV}
Let $[[M_n]]$ and $\langle \langle M_n \rangle \rangle$ be the
matrix-valued quadratic-covariation processes defined above; let
$Ce$ be the matrix-valued deterministic limit defined above. Then
condition {\em \eqn{ek2}} is equivalent to \beql{abc1} [[M_n]]
\Rightarrow C e \qinq D^{k^2} ~, \eeq while condition {\em
\eqn{ek6}} is equivalent to \beql{abc2} \langle \langle M_n
\rangle \rangle \Rightarrow C e \qinq D^{k^2} ~, \eeq
\end{lemma}

Lemma \ref{lemModQCV} is important, not only for general
understanding, but because it plays an important role in the proof
of Theorem \ref{thMart}.  (See the discussion after \eqn{ek101}.)

\paragraph{Proof of Lemma \ref{lemModQCV}.}  We exploit the fact that
the ordinary quadratic variation processes are nondecreasing, but
we need to do more to treat the quadratic covariation processes
when $i \not= j$. Since $[M_{n,i}, M_{n,i}]$ and $\langle M_{n,i},
M_{n,i}\rangle$ for each $i$ are nondecreasing and $\{c_{i,i} t: t
\ge 0\}$ is continuous, we can apply \S VI.2.b and Theorem VI.3.37
of JS to get convergence of these one-dimensional
quadratic-variation processes in $D$ for each $i$ from the
corresponding limits in $\RR$ for each $t$. Before applying
Theorem VI.3.37 of JS, we note that we have convergence of the
finite-dimensional distributions, because we can apply Theorem
11.4.5 of Whitt \cite{W02}.
 We then use the representations
\begin{eqnarray}
2 [M_{n,i}, M_{n,j}] & = & [M_{n,i} + M_{n,j}, M_{n,i}+ M_{n,j}] - [M_{n,i}, M_{n,i}] - [M_{n,j}, M_{n,j}] \nonumber \\
2 \langle M_{n,i}, M_{n,j}\rangle & = & \langle M_{n,i} + M_{n,j},
M_{n,i} + M_{n,j}\rangle
 - \langle M_{n,i}, M_{n,i}\rangle - \langle M_{n,j}, M_{n,j}\rangle ~, \nonumber
\end{eqnarray}
e.g., see \S 1.8 of Liptser and Shiryaev \cite{LS89}.  First, we obtain
the limits \beq [M_{n,i} + M_{n,j}, M_{n,i}+ M_{n,j}] (t)
\Rightarrow 2 c_{i,j} t + c_{i,i} t + c_{j,j} t \qinq \RR \eeqno
and \beq \langle M_{n,i} + M_{n,j}, M_{n,i} + M_{n,j}\rangle  (t)
\Rightarrow 2 c_{i,j} t + c_{i,i} t + c_{j,j} t \qinq \RR \eeqno
for each $t$ and $i \not=j$ from conditions \eqn{ek2} and
\eqn{ek6}. The limits in \eqn{ek2} and \eqn{ek6} for the
components extend to vectors and then we can apply the continuous
mapping theorem with addition. Since $[M_{n,i} + M_{n,j}, M_{n,i}+
M_{n,j}]$ and $\langle M_{n,i} + M_{n,j}, M_{n,i} +
M_{n,j}\rangle$ are both nondecreasing processes, we can repeat
the argument above for $[M_{n,i}, M_{n,i}]$ and $\langle M_{n,i},
M_{n,i}\rangle$ to get \beq [M_{n,i} + M_{n,j}, M_{n,i}+ M_{n,j}]
\Rightarrow (2 c_{i,j} + c_{i,i}  + c_{j,j}) e \qinq D \eeqno and
\beq \langle M_{n,i} + M_{n,j}, M_{n,i} + M_{n,j}\rangle
\Rightarrow (2 c_{i,j} + c_{i,i}  + c_{j,j}) e \qinq D \eeqno We
then get the corresponding limits in $D^3$ for the vector
processes, and apply the continuous mapping theorem with addition
again to get \beq [M_{n,i}, M_{n,j}] \Rightarrow c_{i,j} e \qandq
\langle M_{n,i}, M_{n,j}\rangle \Rightarrow c_{i,j} e \qinq D
\eeqno for all $i$ and $j$.  Since these limits extend to vectors,
we finally have derived the claimed limits in \eqn{abc1} and
\eqn{abc2}.~~~\bsq

\textbf{Outline of the Proof}

The broad outline of the proof of Theorem \ref{thMart} is
standard.  As in both EK and JS, the proof is an application of
Corollary \ref{corProhorov}:  We first show that the sequence
$\{M_n: n \ge 1\}$ is tight in $D^k$, which implies relative
compactness by Theorem \ref{thmProhorov}.
Having established relative compactness, we show that the limit of
any convergent subsequence must be $k$-dimensional Brownian motion
with covariance matrix $C$. That is, there is a tightness step and
there is a characterization step. Both EK and JS in their
introductory remarks, emphasize the importance of the
characterization step.

Before going into details, we indicate where the key results are
in JS. Case (i) in Theorem \ref{thMart} is covered by Theorem
VIII.3.12 on p 432 of JS.  Condition 3.14 on top of p 433 in JS is
implied by condition \eqn{ek1}.  Condition (b.ii) there in JS is
condition \eqn{ek2}. The counterexample in Remark 3.19 on p 434 of
JS shows that weakening condition 3.14 there can cause problems.

Case (ii) in Theorem \ref{thMart} is covered by Theorem VIII.3.22
on p 435 of JS.  Condition 3.23 on p 435 of JS is implied by
condition \eqn{ek5}.  There $\nu$ is the predictable random
measure, which is part of the characteristics of a semimartingale,
as defined in \S II.2 of JS and $*$ denotes the operator in 1.5 on
p 66 of JS constructing the associated integral process.  (We will
not use the notions $\nu$ and $*$ here.)

\section{Tightness}\label{secTightAll}

\subsection{Basic Properties}\label{secTight}
\hsp We work in the setting of a complete separable metric space
(CSMS), also known as a Polish space; see \S\S 13 and 19 of
Billingsley \cite{B99}, \S\S 3.8-3.10 of EK \cite{EK86}
and \S\S 11.1 and 11.2 of Whitt \cite{W02}. (The space $D^k \equiv
D([0, \infty), \RR)^k$ is made a CSMS in a standard way and the
space of probability measures on $D^k$ becomes a CSMS as well.)
Key concepts are:  closed, compact, tight, relatively compact and
sequentially compact. We assume knowledge of metric spaces and
compactness in metric spaces.

\begin{definition} {\em $($tightness$)$}\label{defTight}
A set $A$ of probability measures on a metric space $S$ is
\textbf{tight} if, for all $\epsilon > 0$, there exists a compact
subset $K$ of $S$ such that \beq P(K) > 1 - \epsilon \qforallq P
\in A ~. \eeqno A set of random elements of the metric space $S$
is tight if the associated set of their probability laws on $S$ is
tight. Consequently, a sequence $\{X_n: n \ge 1\}$ of random
elements of the metric space $S$ is tight if, for all $\epsilon >
0$, there exists a compact subset $K$ of $S$ such that \beq P(X_n
\in K) > 1 - \epsilon \qforallq n \ge 1 ~. \eeqno
\end{definition}

Since a continuous image of a compact subset is compact, we have
the following lemma.

\begin{lemma}{\em $($continuous functions of random elements$)$}\label{lemTightCont}
Suppose that $\{X_n: n \ge 1\}$ is a tight sequence of random
elements of the metric space $S$. If $f: S \ra S'$ is a continuous
function mapping the metric space $S$ into another metric space
$S'$, then $\{f(X_n): n \ge 1\}$ is a tight sequence of random
elements of the metric space $S'$.
\end{lemma}

\paragraph{Proof.}  As before, let $\circ$ be used for composition:  $(f \circ g) (x) \equiv f(g(x))$.
 For any function $f: S \ra S'$ and any subset $A$ of $S$, $A \subseteq f^{-1} \circ f (A)$.
Let $\epsilon > 0$ be given. Since $\{X_n: n \ge 1\}$ is a tight
sequence of random elements of the metric space $S$, there exists
a compact subset $K$ of $S$ such that \beq P(X_n \in K) > 1 -
\epsilon \qforallq n \ge 1 ~. \eeqno Then $f(K)$ will serve as the
desired compact set in $S'$, because \beq P(f(X_n) \in f(K)) =
P(X_n \in (f^{-1} \circ f)(K)) \ge  P(X_n \in K) > 1 - \epsilon
\eeqno for all $n \ge 1$.~~~\bsq

We next observe that on products of separable metric spaces
tightness is characterized by tightness of the components; see \S
11.4 of \cite{W02}.

\begin{lemma}{\em $($tightness on product spaces$)$}\label{lemTightProd}
Suppose that $\{(X_{n,1}, \ldots , X_{n,k}) : n \ge 1\}$ is a
sequence of random elements of the product space $S_1 \times
\cdots \times S_k$, where each coordinate space $S_i$ is a
separable metric space. The sequence $\{(X_{n,1}, \ldots ,
X_{n,k}) : n \ge 1\}$ is tight if and only if the sequence
$\{X_{n,i} : n \ge 1\}$ is tight for each $i$, $1 \le i \le k$.
\end{lemma}

\paragraph{Proof.}
The implication from the random vector to the components follows
from Lemma \ref{lemTightCont} because the component $X_{n,i}$ is
the image of the projection map $\pi_i: S_1 \times \cdots \times
S_k \ra S_i$ taking $(x_1, \ldots , x_k)$ into $x_i$, and the
projection map is continuous.  Going the other way, we use the
fact that \beq A_1 \times \cdots \times A_k = \bigcap_{i = 1}^{k}
\pi^{-1}_i (A_i) = \bigcap_{i = 1}^{k} \pi^{-1}_i \circ \pi_i (A_1
\times \cdots \times A_k) \eeqno for all subsets $A_i \subseteq
S_i$.  Thus, for each $i$ and any $\epsilon > 0$, we can choose
$K_i$ such that $P(X_{n,i} \notin K_i) < \epsilon/k$ for all $n
\ge 1$.  We then let $K_1 \times \cdots \times K_k$ be the desired
compact for the random vector.  We have
\begin{eqnarray}
P\left((X_{n,1}, \ldots , X_{n,k}) \notin K_1 \times \cdots \times K_k\right) & = & P\left(\bigcup_{i=1}^{k} \{X_{n,i} \notin K_i\} \right) \nonumber \\
& \le & \sum_{i=1}^{k} P\left( X_{n,i} \notin K_i \right) \le
\epsilon ~.~~~\bsq \nonumber
\end{eqnarray}




Tightness goes a long way toward establishing convergence because
of Prohorov's theorem.  It involves the notions of sequential
compactness and relative compactness.

\begin{definition} {\em $($relative compactness and sequential compactness$)$}\label{defRelCompact}
A subset $A$ of a metric space $S$ is \textbf{relatively compact}
if every sequence $\{x_n: n \ge 1\}$ from $A$ has a subsequence
that converges to a limit in $S$ $($which necessarily belongs to
the closure $\bar{A}$ of $A)$. A subset of $S$ is
\textbf{sequentially compact} if it is closed and relatively
compact.
\end{definition}

We rely on the following basic result about compactness on metric
spaces.

\begin{lemma} {\em $($compactness coincides with sequential compactness on metric spaces$)$}\label{defRelCompact2}
A subset $A$ of a metric space $S$ is compact if and only if it is
sequentially compact.
\end{lemma}

We can now state Prohorov's theorem; see \S 11.6 of \cite{W02}.
It relates compactness of sets of measures to compact subsets of
the underlying sample space $S$ on which the probability measures
are defined.

\begin{theorem}{\em $($Prohorov's theorem$)$}\label{thmProhorov}
A subset of probability measures on a CSMS is tight if and only if
it is relatively compact.
\end{theorem}

We have the following elementary corollaries:

\begin{coro}{\em $($convergence implies tightness$)$}\label{corConvTight}
If $X_n \Rightarrow X$ as $n \ra \infty$ for random elements of a
CSMS, then the sequence $\{X_n: n \ge 1\}$ is tight.
\end{coro}

\begin{coro}{\em $($individual probability measures$)$}\label{corIndiv}
Every individual probability measure on a CSMS is tight.
\end{coro}

As a consequence of Prohorov's Theorem, we have the following
method for establishing convergence of random elements:

\begin{coro}{\em $($convergence in distribution via tightness$)$}\label{corProhorov}
Let $\{X_n: n \ge 1\}$ be a sequence of random elements of a CSMS
$S$.  We have \beq X_n \Rightarrow X \qinq S \qasq n \ra \infty
\eeqno if and only if (i) the sequence $\{X_n: n \ge 1\}$ is tight
and (ii) the limit of every convergent subsequence of $\{X_n: n
\ge 1\}$ is the same fixed random element $X$ (has a common
probability law).
\end{coro}

In other words, once we have established tightness, it only
remains to show that the limits of all converging subsequences
must be the same. With tightness, we only need to uniquely
determine the limit.  When proving Donsker's theorem, it is
natural to uniquely determine the limit through the
finite-dimensional distributions.  Convergence of all the
finite-dimensional distributions is not enough to imply
convergence on $D$, but it does uniquely determine the
distribution of the limit; see pp 20 and 121 of Billingsley \cite{B68}
and Example 11.6.1 in Whitt \cite{W02}.

We will apply this approach to prove the martingale FCLT in this paper.
In the martingale setting it is natural
to use the martingale characterization of Brownian motion,
originally established by L\'{e}vy \cite{L48} and proved by Ito's
formula by Kunita and Watanabe \cite{KW67}; see p. 156 of Karatzas and
Shreve \cite{KS88}, and various extensions, such as to continuous
processes with independent Gaussian increments, as in Theorem 1.1
on p. 338 of EK \cite{EK86}. A thorough study of
martingale characterizations appears in Chapter 4 of Liptser and
Shiryayev \cite{LS89} and in Chapters VIII and IX of JS \cite{JS87}.

\subsection{Stochastic Boundedness}\label{secSB}
\hsp
We now discuss stochastic boundedness because it plays a role in the tightness
criteria in the next section.
We start by defining stochastic boundedness and relating it to tightness.
We then discuss situations in which stochastic boundedness is preserved.
Afterwards, we give conditions for a sequence of martingales to be stochastically bounded in $D$
involving the stochastic boundedness of appropriate sequences of $\RR$-valued random variables.

\subsubsection{Connection to Tightness}

 For random elements of $\RR$ and $\RR^k$, stochastic boundedness and tightness
 are equivalent, but tightness is stronger than stochastic boundedness for
random elements of the functions spaces $C$ and $D$ (and the associated product spaces $C^k$ and $D^k$).

\begin{definition} {\em $($stochastic boundedness for random vectors$)$}\label{defSB}
A sequence $\{X_n: n \ge 1\}$ of random vectors taking values in $\RR^k$ is \textbf{stochastically bounded (SB)} if the sequence is tight, as defined in
Definition {\em \ref{defTight}}.
\end{definition}

The notions of tightness and stochastic boundedness thus agree for random elements of $\RR^k$, but these
notions differ for stochastic processes.
For a function $x \in D^k \equiv D([0, \infty), \RR)^k$, let
\beq
\| x\|_T \equiv \sup_{0 \le t \le T}{\{ | x (t) |\}} ~,
\eeqno
where $|b|$ is a norm of $b \equiv (b_1, b_2, \ldots , b_k)$ in $\RR^k$ inducing the Euclidean topology, such as the maximum norm:
$|b| \equiv \max{\{|b_1|, |b_2|, \ldots , |b_k|\}}$. (Recall that all norms on Euclidean space $\RR^k$ are equivalent.)

\begin{definition} {\em $($stochastic boundedness for random elements of $D^k)$}\label{defSBD}
A sequence $\{X_n: n \ge 1\}$ of random elements of $D^k$ is \textbf{stochastically bounded in $D^k$} if the sequence
of real-valued random variables $\{\| X_n\|_T: n \ge 1\}$
is stochastically bounded in $\RR$ for each $T > 0$, using Definition {\em \ref{defSB}}.
\end{definition}

For random elements of $D^k$, tightness is a strictly stronger concept than stochastic boundedness.
Tightness of $\{X_n\}$ in $D^k$
implies stochastic boundedness, but not conversely; see \S 15 of Billingsely \cite{B68}.

\subsubsection{Preservation}

We have the following analog of Lemma \ref{lemTightProd}, which characterizes tightness for sequences of random vectors in terms
of tightness of the associated sequences of components.

\begin{lemma}{\em $($stochastic boundedness on $D^k$ via components$)$}\label{lemSBprod}
A sequence
$$\{(X_{n,1}, \ldots , X_{n,k}) : n \ge 1\} \qinq D^k \equiv D \times \cdots \times D$$
is stochastically bounded in $D^k$ if and only if the sequence $\{X_{n,i} : n \ge 1\}$
is stochastically bounded in $D \equiv D^1$ for each $i$,
$1 \le i \le k$.
\end{lemma}

\paragraph{Proof.}  Assume that we are using the maximum norm on product spaces.
We can apply Lemma \ref{lemTightProd} after noticing that
\beq
\| (x_1, \ldots , x_k)\|_T = \max{\{ \| x_i\|_T: 1 \le i \le k\}}
\eeqno
for each element $(x_1, \ldots , x_k)$ of $D^k$.
Since other norms are equivalent, the result applies more generally.~~~\bsq

\begin{lemma}{\em $($stochastic boundedness in $D^k$ for sums$)$}\label{lemSBsums}
Suppose that
\beq
Y_n (t) \equiv X_{n,1} (t) + \cdots + X_{n,k} (t), \quad t \ge 0,
\eeqno
for each $n \ge 1$, where
$\{(X_{n,1}, \ldots , X_{n,k}) : n \ge 1\}$ is a sequence of random elements of the product space $D^k \equiv D \times \cdots \times D$.
If $\{X_{n,i} : n \ge 1\}$
is stochastically bounded in $D$ for each $i$,
$1 \le i \le k$, then the sequence $\{Y_n: n \ge 1\}$
is stochastically bounded in $D$.
\end{lemma}

Note that the converse is not true:  We could have $k=2$ with $X_{n,2} (t) = - X_{n,1} (t)$ for all $n$ and $t$.
In that case we have $Y_n (t) = 0$ for all $X_{n,1} (t)$.

\subsubsection{Stochastic Boundedness for Martingales}

We now provide ways to get stochastic boundedness for sequences of martingales in $D$
from associated sequences of random variables.  Our first result exploits the classical submartingale-maximum inequality;
e.g., see p. 13 of Karatzas and Shreve \cite{KS88}.  We say that a function $f: \RR \ra \RR$ is {\em even} if
$f(-x) = f(x)$ for all $x \in \RR$.

\begin{lemma}{\em $($SB from the maximum inequality$)$}\label{lemSBmax}
Suppose that, for each $n \ge 1$, $M_n \equiv \{M_n (t): t \ge 0\}$ is a martingale $($with respect to a specified filtration$)$
with sample paths in $D$.
Also suppose that, for each $T > 0$, there exists an even nonnegative convex function $f: \RR \ra \RR$
with first derivative $f' (t) > 0$ for $t > 0$ $($e.g., $f(t) \equiv t^2)$,
there exists a positive constant $K \equiv K(T, f)$, and there exists
an integer $n_0 \equiv n_0 (T,f, K)$,
 such that
\beq
E[f(M_n(T))] \le K \qforallq n \ge n_0 ~.
\eeqno
Then the sequence of stochastic processes $\{M_n: n \ge 1\}$ is stochastically bounded in $D$.
 \end{lemma}

 \paragraph{Proof.}  Since any set of finitely many random elements of $D$ is automatically tight,
 Theorem 1.3 of Billingsley \cite{B99}, it suffices to consider $n \ge n_0$.
 Since $f$ is continuous and $f'(t) >0 $ for $t > 0$, $t > c$ if and only if
 $f(t) > f(c)$ for $t > 0$.  Since $f$ is even,
 \beq
 E[f(M_n(t))] = E[f(| M_n(t)| )] \le E[f(|M_n(T)|)] = E[f(M_n(T))]\le K
 \eeqno
 for all $t$, $0 \le t \le T$.
 Since these moments are finite and $f$ is convex,
 the stochastic process $\{f(M_n(t)) : 0 \le t \le T\}$ is a submartingale for each $n \ge 1$, so that
 we can apply the submartingale-maximum inequality to get
 \beq
 P(\| M_n \|_T > c) =  P(\| f \circ M_n \|_T > f(c)) \le \frac{E[f(M_n(T))]}{f(c)} \le \frac{K}{f(c)}
 \eeqno
 for all $n \ge n_0$.
 Since $f(c) \ra \infty$ as $c \ra \infty$, we have the desired conclusion.~~~\bsq

We now establish another sufficient condition for stochastic boundedness of
square-integrable martingales by applying the Lenglart-Rebolledo inequality;
see p. 66 of Liptser and Shiryayev \cite{LS89} or p. 30 of Karatzas and Shreve \cite{KS88}.

\begin{lemma}{\em $($Lenglart-Rebolledo inequality$)$}\label{lemLenglart}
Suppose that $M \equiv \{M (t): t \ge 0\}$ is a square-integrable martingale $($with respect to a specified filtration$)$
with predictable quadratic variation $\langle M \rangle \equiv \{\langle M \rangle (t): t \ge 0\}$, i.e.,
such that $M^2 - \langle M \rangle \equiv \{ M (t)^2 - \langle M \rangle (t): t \ge 0\}$ is a martingale
by the Doob-Meyer decomposition.  Then, for all $c>0$ and $d>0$,
\beql{lenglart}
P \left( \sup_{0 \le t \le T}{\{ |M(t)|\}} > c \right) \le \frac{d}{c^2} + P\left(\langle M \rangle (T) > d \right) ~.
\eeq
 \end{lemma}

As a consequence we have the following criterion for stochastic boundedness of a sequence of square-integrable martingales.

\begin{lemma}{\em $($SB criterion for square-integrable martingales$)$}\label{lemSBLenglart}
Suppose that, for each $n \ge 1$, $M_n \equiv \{M_n (t): t \ge 0\}$ is a square-integrable martingale $($with respect to a specified filtration$)$
with predictable quadratic variation $\langle M_n \rangle \equiv \{\langle M_n \rangle (t): t \ge 0\}$, i.e.,
such that $M_n^2 - \langle M_n \rangle \equiv \{ M_n (t)^2 - \langle M_n \rangle (t): t \ge 0\}$ is a martingale
by the Doob-Meyer decomposition.  If the sequence of random variables $\{\langle M_n \rangle (T): n \ge 1\}$
is stochastically bounded in $\RR$ for each $T > 0$, then the sequence of stochastic processes $\{M_n: n \ge 1\}$ is stochastically bounded
in $D$.
 \end{lemma}

\paragraph{Proof.}
For $\epsilon > 0$ given, apply the assumed stochastic boundedness of
the sequence $\{\langle M_n \rangle (T): n \ge 1\}$ to obtain a constant $d$ such that
\beq
P\left(\langle M_n \rangle (T) > d \right) < \epsilon/2 \qforallq n \ge 1 ~.
\eeqno
Then for that determined $d$, choose $c$ such that $d/c^2 < \epsilon/2$.  By the
Lenglart-Rebolledo inequality \eqn{lenglart}, these two inequalities imply that
\beql{SBlenglart}
P \left( \sup_{0 \le t \le T}{\{ |M_n(t)|\}} > c \right) < \epsilon ~.~~~\bsq
\eeq

\subsection{Tightness Criteria}\label{secSimp}
\hsp
The standard characterization for tightness of a sequence of
stochastic processes in $D$, originally developed by Skorohod
\cite{S56} and presented in Billingsley \cite{B68}, involves suprema.
Since then, more elementary one-dimensional criteria have been
developed; see Billingsley \cite{B74}, Kurtz \cite{K75}, Aldous \cite{A78, A89},
Jacod et al. \cite{JMM83}, \S\S 3.8 and 3.9 of EK \cite{EK86}
and \S 16 of Billingsley \cite{B99}.  Since these
simplifications help in proving the martingale FCLT, we will
present some of the results here.

\subsubsection{Criteria Involving a Modulus of Continuity}\label{secClassic}

We start by presenting the classical characterization of
tightness, as in Theorems 13.2 and 16.8 of Billingsley \cite{B99}. For
that purpose we define functions $w$ and $w'$ that can serve as a
modulus of continuity. For any $x \in D$ and subset $A$ of $[0,
\infty)$, let \beql{mod1} w(x, A) \equiv \sup{\{|x(t_1) - x(t_2)|:
t_1, t_2 \in A\}} ~. \eeq For any $x \in D$, $T > 0$ and $\delta >
0$, let \beql{mod1a} w(x,\delta, T ) \equiv \sup{\{ w(x, [t_1,
t_2]): 0 \le t_1 < t_2 \le (t_1 + \delta) \wedge T\}} \eeq and
\beql{mod2} w'(x,\delta, T ) \equiv \inf_{\{ t_i\}}  \max_{1 \le i
\le k}{\{w(x, [t_{i-1}, t_i))\}} ~, \eeq where the infimum in
\eqn{mod2} is over all $k$ and all subsets of $[0,T]$ of size
$k+1$ such that \beq 0 = t_0 < t_1 < \cdots < t_k = T \quad
\mbox{with} \quad t_i - t_{i-1} > \delta \qforq 1 \le i \le k-1 ~.
\eeqno (We do not require that $t_k - t_{k-1} > \delta$.)

The following is a variant of the classical characterization of
tightness; see Theorem 16.8 of Billingsley \cite{B99}.
\begin{theorem}{\em $($classical characterization of tightness$)$}\label{thClassic}
A sequence of stochastic processes $\{X_n: n \ge 1\}$ in $D$ is
tight if and only if

\vspace{0.03in} $(i)$ The sequence $\{X_n: n \ge 1\}$ is
stochastically bounded in $D$

and

\vspace{0.03in} $(ii)$ for each $T > 0$ and $\epsilon > 0$, \beq
\lim_{\delta \downarrow 0} \limsup_{n \ra \infty} P(w'(X_n,
\delta, T) > \epsilon) = 0 ~. \eeqno
 \end{theorem}

If the modulus $w (X_n, \delta, T)$ is substituted for $w'(X_n,
\delta, T)$ in Condition (ii)
 of Theorem \ref{thClassic}, then the sequence $\{X_n\}$ is said to be \textbf{C-tight}, because then
 the sequence is again tight but the limit of any convergent subsequence must have continuous sample paths; see
 Theorem 15.5 of Billingsley \cite{B68}.  With this modified condition (ii), condition (i) is implied
 by having the sequence $\{X_n (0)\}$ be tight in $\RR$.

 Conditions (i) and (ii) in Theorem \ref{thClassic} are both somewhat hard to verify because they involve suprema.  We have shown how the
 stochastic-boundedness condition (i) can be made one-dimensional for martingales via Lemmas \ref{lemSBmax}
 and \ref{lemSBLenglart}.  The following is another result, which exploits the modulus condition (ii) in Theorem \ref{thClassic};
 see p. 175 of Billingsley \cite{B99}.  To state it, let $J$ be the maximum-jump function defined in \eqn{mod3}.

\begin{lemma}{\em $($substitutes for stochastic boundedness$)$}\label{lemSub1}
In the presence of the modulus condition $(ii)$ in Theorem {\em
\ref{thClassic}}, each of the following is equivalent to the
stochastic-boundedness condition $(i)$ in Theorem {\em
\ref{thClassic}}:

\vspace{0.03in} $(i)$ The sequence $\{X_n (t): n \ge 1\}$ is
stochastically bounded in $\RR$ for each $t$ in a dense subset of
$[0,\infty)$.

\vspace{0.03in} $(ii)$ The sequence $\{X_n (0): n \ge 1\}$ is
stochastically bounded in $\RR$ and, for each $T > 0$, the
sequence $\{J(X_n,T): n \ge 1\}$ is stochastically bounded in
$\RR$.
 \end{lemma}

We also mention the role of the maximum-jump function $J$ in
\eqn{mod3}
 in characterizing continuous limits; see Theorem 13.4 of Billingsley \cite{B99}:

\begin{lemma}{\em $($identifying continuous limits$)$}\label{lemContLim}
Suppose that $X_n \Rightarrow X$ in $D$. Then $P(X \in C) = 1$ if
and only if $J(X_n,T) \Rightarrow 0$ in $\RR$ for each $T>0$.
 \end{lemma}

\subsubsection{Simplifying the Modulus Condition}\label{secModulusSimp}

Simplifying the modulus condition (ii) in Theorem \ref{thClassic}
is even of greater interest. We first present results from EK.
  Conditions (i) and (iii) below are of
particular interest because they shows that, for each $n$ and $t$,
we need only consider the process at the two single time points
$t+u$ and $t-v$ for $u >0$ and $v>0$.  As we will see in Lemma
\ref{lemSimple} below, we can obtain a useful sufficient condition
involving only the single time point $t+u$ (ignoring $t-v$).

For the results in this section, we assume that the strong
stochastic-bounded\-ness condition (i) in Theorem \ref{thClassic} is
in force.  For some of the alternatives to the modulus condition
(ii) in Theorem \ref{thClassic} it is also possible to simplify
the stochastic-boundedness condition, as in Lemma \ref{lemSub1},
but we do not carefully examine that issue.

\begin{theorem}{\em $($substitutes for the modulus condition$)$}\label{thSub2}
In the presence of the stochastic-boundedness condition $(i)$ in
Theorem {\em \ref{thClassic}}, \textbf{each of the following is
equivalent to the modulus condition} $(ii)$:

\vspace{0.03in} $(i)$ For each $T > 0$, there exists a constant
$\beta > 0$ and a family of nonnegative random variables $\{Z_n
(\delta, T): n \ge 1, \delta > 0\}$ such that
\begin{eqnarray}\label{mod6}
&& E\left[ (1 \wedge |X_n (t+ u) - X_n (t)|)^{\beta}\left| \right. \sF_{n, t}\right] \left((1 \wedge|X_n (t) - X_n (t - v)|)^{\beta}\right) \nonumber \\
&& \quad \quad \le E\left[Z_n (\delta, T)| \sF_{n,t}\right] \quad
\mbox{w.p.1}
\end{eqnarray}
for $0 \le t \le T$, $0 \le u \le \delta$ and $0 \le v \le
t\wedge\delta$, where $\sF_{n,t}$ is the $\sigma$-field in the
internal filtration of $\{X_n (t): t \ge 0\}$, \beql{mod7}
\lim_{\delta \downarrow 0} \limsup_{n \ra \infty} E[Z_n (\delta,
T)] = 0 ~. \eeq and \beql{mod7a} \lim_{\delta \downarrow 0}
\limsup_{n \ra \infty} E[(1 \wedge |X_n (\delta) - X_n
(0)|)^{\beta}] = 0 ~. \eeq

\vspace{0.03in} $(ii)$ The sequence of stochastic processes
$\{\{f(X_n (t)): t \ge 0\}: n \ge 1\}$ is tight in $D$ for each
function $f: \RR \ra \RR$ in a dense subset of all bounded
continuous functions in the topology of uniform convergence over
bounded intervals.

\vspace{0.03in} $(iii)$ For each function $f: \RR \ra \RR$ in a
dense subset of all bounded continuous functions in the topology
of uniform convergence over bounded intervals, and $T > 0$, there
exists a constant $\beta > 0$ and a family of nonnegative random
variables $\{Z_n (\delta, f, T): n \ge 1, \delta > 0\}$ such that
\begin{eqnarray}\label{mod7b}
&& E\left[ |f(X_n (t+ u)) - f(X_n (t))|^{\beta}\left| \right. \sF_{n,f,t}\right] \left(|f(X_n (t)) - f(X_n (t - v))|^{\beta}\right) \nonumber \\
&& \quad \quad \le E\left[Z_n (\delta, f, T)| \sF_{n,f,t}\right]
\quad \mbox{w.p.1}
\end{eqnarray}
 for $0 \le t \le T$, $0 \le u \le \delta$ and $0 \le v \le t\wedge\delta$, where $\sF_{n,f,t}$ is the $\sigma$-field
in the internal filtration of $\{f(X_n (t)): t \ge 0\}$, and
\beql{mod7c} \lim_{\delta \downarrow 0} \limsup_{n \ra \infty}
E[Z_n (\delta, f, T)] = 0 ~. \eeq and \beql{mod7d} \lim_{\delta
\downarrow 0} \limsup_{n \ra \infty} E[  |f(X_n (\delta)) - f(X_n
(0))|^{\beta}] = 0 ~. \eeq

 \end{theorem}

 \paragraph{Proof.}  Condition (i) is condition (b) in Theorem 3.8.6 of EK
 for the case of real-valued stochastic processes.
 Condition (ii) is Theorem 3.9.1 of EK.
 Condition (iii) is a modification of condition (b) in Theorem 3.8.6 of EK, exploiting condition (ii) and
the fact that our stochastic processes are real valued.  Since the
functions $f$ are bounded in (iii), there is no need to replace
the usual metric $m (c,d) \equiv |c - d|$ with the bounded metric
$q \equiv m\wedge 1$.~~~\bsq

In applications of Theorem \ref{thSub2}, it is advantageous to
replace Conditions (i) and (iii) in Theorem \ref{thSub2} with
stronger sufficient conditions, which only involves the
conditional expectation on the left.  By doing so, we need to
consider only conditional distributions of the processes at a
single time point $t + u$ with $u>0$, given the relevant history
$\sF_{n,t}$ up to time time $t$.  We need to do an estimate for
all $t$ and $n$, but given specific values of $t$ and $n$, we need
to consider only one future time point $t+u$ for $u > 0$.
\textbf{These simplified conditions are sufficient, but not
necessary, for $D$-tightness}.  On the other hand, they are not
sufficient for $C$-tightness; see Remark \ref{rmSimple} below.

In addition, for condition (iii) in Theorem \ref{thSub2}, it
suffices to specify a specific dense family of functions. As in
the Ethier-Kurtz proof of the martingale FCLT (their Theorem
7.1.4), we introduce smoother functions in order to exploit Taylor
series expansions. Indeed, we present the versions of Theorem
\ref{thSub2} actually applied by EK in their
proof of their Theorem 7.1.4.

\begin{lemma}{\em $($simple sufficient criterion for tightness$)$}\label{lemSimple}
The following are sufficent, but not necessary, conditions for a
sequence of stochastic processes $\{X_n: n \ge 1\}$ in $D$ to be
tight:

\vspace{0.03in} $(i)$ The sequence $\{X_n: n \ge 1\}$ is
stochastically bounded in $D$

\vspace{0.03in} and either

\vspace{0.03in} $(ii.a)$  For each $n \ge 1$, the stochastic
process $X_n$ is adapted to a filtration \textbf{F}$_n \equiv
\{\sF_{n,t}: t \ge 0\}$. In addition, for each function $f: \RR
\ra \RR$ belonging to $\sC^{\infty}_{c}$ (having compact support
and derivatives of all orders) and $T > 0$, there exists a family
of nonnegative random variables $\{Z_n (\delta, f, T): n \ge 1,
\delta > 0\}$ such that \beql{mod8} \left| E\left[ f(X_n (t+ u)) -
f(X_n (t)) \left| \right. \sF_{n,t} \right] \right|  \le
E\left[Z_n (\delta, f, T)\left| \right. \sF_{n,t}\right] \eeq
w.p.1 for $0 \le t \le T$ and $0 \le u \le \delta$ and \beql{mod9}
\lim_{\delta \downarrow 0} \limsup_{n \ra \infty} E[Z_n (\delta,
f, T)] = 0 ~. \eeq or

\vspace{0.03in} $(ii.b)$  For each $n \ge 1$, the stochastic
process $X_n$ is adapted to a filtration \textbf{F}$_n \equiv
\{\sF_{n,t}: t \ge 0\}$. In addition, for each $T > 0$, there
exists a family of nonnegative random variables $\{Z_n (\delta,
T): n \ge 1, \delta > 0\}$ such that \beql{mod9a} \left| E\left[
\left(X_n (t+ u) - X_n (t)\right)^2 \left| \right. \sF_{n,t}
\right] \right|  \le E\left[Z_n (\delta, T)\left| \right.
\sF_{n,t}\right] \eeq w.p.1 for $0 \le t \le T$ and $0 \le u \le
\delta$ and \beql{mod9b} \lim_{\delta \downarrow 0} \limsup_{n \ra
\infty} E[Z_n (\delta, T)] = 0 ~. \eeq
\end{lemma}

\paragraph{Proof.}  We apply Theorem \ref{thSub2}.  For condition (ii.a), we apply
Theorem \ref{thSub2} (iii). We have specified a natural class of
functions $f$ that are dense in the set of all continuous bounded
real-valued functions with the topology of uniform convergence
over bounded intervals.  We have changed the filtration.  Theorem
\ref{thSub2} and Theorem 3.8.6 of EK specify
the internal filtration of the process, which here would be the
process $\{f(X_n (t)): t \ge 0\}$, but it is convenient to work
with the more refined filtration associated with $X_n$.  By taking
conditional expected values, conditional on the coarser filtration
generated by $\{f(X_n (t)): t \ge 0\}$, we can deduce the
corresponding inequality conditioned on the coarser filtration.
The conditions here follow from condition \eqn{mod6} by taking
$\beta = 2$, because
\begin{eqnarray}\label{mod10}
&& E[ |f(X_n (t+ u)) - f(X_n (t))|^{2}| \sF_{n,t}]  =   E[ f(X_n (t+ u))^2 - f(X_n (t))^2 | \sF_{n,t}]  \nonumber \\
&& \quad \quad - 2 f(X_n (t))E[ f(X_n (t+ u)) - f(X_n (t))|
\sF_{n,t}] ~, \nonumber
\end{eqnarray}
(see (1.35) on p. 343 of EK), so that
\beq E[ |f(X_n (t+
u)) - f(X_n (t))|^{2}| \sF_{n,t}] \le E[Z_n(\delta, f^2, T)|
\sF_{n,t}] + 2 \|f\| E[Z_n(\delta, f, T)| \sF_{n,t}]
\eeqno
under
condition \eqn{mod8}.  Note that $f^2$ belongs to our class of
functions for every function $f$ in the class. Finally, note that
under \eqn{mod8} \beql{mod11a} E[ |f(X_n (\delta)) - f(X_n
(0))|^{2}] \le E[Z_n (\delta, f, T)] ~, \eeq so that the limit
\eqn{mod7d} holds by \eqn{mod9}.

Now consider condition (ii.b).  Note this is a direct consequence
of Theorem \ref{thSub2} (i) using $\beta = 2$. Condition
\eqn{mod10} is made stronger than \eqn{mod6} by removing the
$\wedge 1$. As in \eqn{mod11a}, \beql{mod11b} E[ |X_n (\delta) -
X_n (0)|^{2} \le E[Z_n (\delta, T)] ~, \eeq so that the limit
\eqn{mod7a} holds by \eqn{mod9b}.
  ~~~\bsq

We have mentioned that Theorem \ref{thSub2} and Lemma
\ref{lemSimple} come from EK \cite{EK86}. As indicated
there, these in turn come from Kurtz \cite{K75}.  However, Lemma
\ref{lemSimple} is closely related to a sufficient condition for
tightness proposed by Billingsley \cite{B74}. Like Lemma
\ref{lemSimple}, this sufficient condition is especially
convenient for treating Markov processes. We state a version of
that here. For that purpose, let $\alpha_n (\lambda, \epsilon,
\delta, T)$ be a number such that \beql{B74a} P( |X_n (u) - X_n
(t_m)| > \epsilon | X_n (t_1), \cdots , X_n (t_m)) \le \alpha_n
(\lambda, \epsilon, \delta, T) \eeq holds with probability $1$ on
the set $\{\max_i \{|X_n (t_i)| \le \lambda \}$ for all $m$ and
all $m$ time points $t_i$ with \beql{B74b} 0 \le t_1 \le \cdots
\le t_{m-1} \le t_m < u \le T \qandq u - t_m \le \delta ~. \eeq
The following is a variant of Theorem 1 of Billingsley \cite{B74}.
(He has the stochastic-boundedness condition (i) below replaced by
a weaker condition.)

\begin{lemma}{\em $($another simple sufficient condition for tightness$)$}\label{lemSimple2}
Alternative sufficent, but not necessary, conditions for a
sequence of stochastic processes $\{X_n: n \ge 1\}$ in $D$ to be
tight are the following:

\vspace{0.03in} $(i)$ The sequence $\{X_n: n \ge 1\}$ is
stochastically bounded in $D$

\vspace{0.03in} and

\vspace{0.03in} $(ii)$ For each $T > 0$, $\epsilon > 0$ and
$\lambda < \infty$, \beql{B74c} \lim_{\delta \downarrow 0}
\limsup_{n \ra \infty} \alpha_n (\lambda, \epsilon, \delta, T) = 0
~, \eeq for $\alpha_n (\lambda, \epsilon, \delta, T)$ defined in
{\em \eqn{B74a}}.
\end{lemma}

\begin{remark}{$($sufficient for $D$-tightness but not $C$-tightness.$)$}\label{rmSimple}
{\em The sufficient conditions in Lemmas \ref{lemSimple} and
\ref{lemSimple2} are sufficient but not necessary for $D$
tightness.  On the other hand, these conditions are {\em not}
sufficient for $C$-tightness.  To substantiate these claims, it
suffices to consider simple examples. The single function $x (t) =
1_{[1,\infty)} (t)$, $t \ge 0$, is necessarily tight in $D$, but
the conditions in Lemmas \ref{lemSimple} and \ref{lemSimple2} are
not satisfied for it. On the other hand, the stochastic process $X
(t) = 1_{[T, \infty)} (t)$, $t \ge 0$, where $T$ is an exponential
random variable with mean $1$, is $D$ tight, but not $C$-tight.
By the lack of memory property of the exponential distribution,
the conditions of Lemmas \ref{lemSimple} and \ref{lemSimple2} are
satisfied for this simple random element of $D$.~~~\bsq }
\end{remark}

\subsubsection{Stopping Times and Quadratic Variations}\label{secStop}

We conclude this section by mentioning alternative criteria for
tightness specifically intended for martingales.  These criteria
involve stopping times and the quadratic-variation processes. The
criteria in terms of stopping times started with Aldous \cite{A78} and
Rebolledo \cite{R80}. Equivalence with the conditions in Theorem
\ref{thSub2} and Lemma \ref{lemSimple} was shown in Theorems 2.7
of Kurtz \cite{K81} and 3.8.6 of EK \cite{EK86}; See also p
176 of Billingsley \cite{B99}.  These criteria are very natural for
proving tightness of martingales, as Aldous \cite{A78, A89} and
Rebolledo \cite{R80} originally showed.  Aldous' \cite{A78} proof of (a
generalization of) the martingale FCLT from McLeish \cite{M74} is
especially nice.

\begin{theorem}{\em $($another substitute for the modulus condition$)$}\label{thSub3}
Suppose that the stochastic-boundedness condition $(i)$ in Theorem
{\em \ref{thClassic}} holds.

\vspace{0.1in} The following is equivalent to the modulus
condition in Theorem {\em \ref{thClassic}}, condition $(ii)$ in
Theorem {\em \ref{thSub2}} and thus to tightness:

\vspace{0.1in} For each $T > 0$, there exists a constant $\beta >
0$ such that \beql{mod301} \lim_{\delta \downarrow 0} \limsup_{n
\ra \infty} C_n (\delta, T) = 0 ~, \eeq where $C_n (\delta, T)$ is
the supremum of \beq E\left[ (1 \wedge |X_n (\tau + u) - X_n
(\tau)|)^{\beta} \left((1 \wedge|X_n (\tau) - X_n (\tau -
v)|)^{\beta}\right)\right] \eeqno over $0 \le u \le \delta$, $\tau
\ge v$ and $\tau \in \sS_{n,T}$, with $\sS_{n,T}$ being the
collection of all finite-valued stopping times with respect to the
internal filtration of $X_n$, bounded by $T$.
\end{theorem}

\begin{theorem}{\em $($another substitute for the sufficient condition$)$}\label{thSub4}
Suppose that the stochastic-boundedness condition $(i)$ in Theorem
{\em \ref{thClassic}} holds.

\vspace{0.1in} Then the following are equivalent sufficient, but
not necessary, conditions for the modulus condition $(ii)$ in
Theorem {\em \ref{thClassic}}, condition $(ii)$ in Theorem {\em
\ref{thSub2}} and thus for tightness:

\vspace{0.03in} $(i)$ For each $T > 0$ and sequence $\{(\tau_n,
\delta_n): n \ge 1\}$, where $\tau_n$ is a finite-valued stopping
time with respect to the internal filtration generated by $X_n$,
with $\tau_n \le T$ and $\delta_n$ is a positive constant with
$\delta_n \downarrow 0$, \beql{mod41} X_n (\tau_n + \delta_n) -
X_n (\tau_n) \Rightarrow 0 \qasq n \ra \infty ~. \eeq

\vspace{0.03in} $(ii)$ For each $\epsilon > 0$, $\eta > 0$ and $T
> 0$, there exist $\delta > 0$, $n_0$ and two sequences
$\{\tau_{n,i}: n \ge n_0\}$, $i = 1, 2$, where $\tau_{n,i}$ is a
stopping time with respect to the internal filtration generated by
$X_n$ such that $0 \le \tau_{n,1} < \tau_{n,2} \le T$ and
\beql{mod42} P \left(X_n (\tau_{n,2}) - X_n (\tau_{n,1}) >
\epsilon, \quad \tau_{n,2} - \tau_{n,1} \le \delta\right) < \eta
\qforallq n \ge n_0 ~. \eeq
 \end{theorem}

Aldous (1978) showed that the two conditions in Lemma
\ref{lemSimple2} imply condition \eqn{mod41}, because \beql{B74d}
P( X_n (\tau_n + \delta_n) - X_n (\tau_n) > \epsilon) \le \alpha_n
(\lambda, \epsilon, \delta_n, T) + P(\|X_n\|_T > \lambda ) ~. \eeq

Related criteria for tightness in terms of quadratic variation
processes are presented in Rebolledo \cite{R80}, Jacod et al. \cite{JMM83}
and \S\S VI.4-5 of JS \cite{JS87}. The following is a
minor variant of Theorem VI.4.13 on p 322 of \cite{JS87}.
Lemma \ref{lemTightProd} shows that $C$-tightness of the
sum $\sum_{i=1}^k \langle X_{n,i}, X_{n,i} \rangle$ is equivalent
to tightness of the components.

\begin{theorem}{\em $($tightness criterion in terms of the angle-bracket process$)$}\label{thAngle}
Suppose that $X_n \equiv (X_{n,1}, \ldots , X_{n,1})$ is a
locally-square-integrable martingale in $D^k$ for each $n\ge 1$,
so that the predictable quadratic covariation processes $\langle
X_{n,i}, X_{n,j} \rangle$ are well defined.  The sequence $\{X_n:
n \ge 1\}$ is $C$-tight if

\vspace{0.03in} $(i)$ The sequence $\{X_{n,i} (0): n \ge 1\}$ is
tight in $\RR$ for each $i$, and

\vspace{0.03in} $(ii)$ the sequence $\langle X_{n,i}, X_{n,i}
\rangle$ is $C$-tight for each $i$.
 \end{theorem}

 The following is a minor variant of Lemma 11 of Rebolledo \cite{R80}.  We again apply Lemma \ref{lemTightProd}
 to treat the vector case.  Since a local martingale with bounded jumps is
 locally square integrable, the hypothesis of the next theorem is not weaker than
 the hypothesis of the previous theorem.

\begin{theorem}{\em $($tightness criterion in terms of the square-bracket processes$)$}\label{thSquare}
Suppose that $X_n \equiv (X_{n,1}, \ldots , X_{n,1})$ is a local
martingale in $D^k$ for each $n\ge 1$.
 The sequence $\{X_n: n \ge 1\}$ is $C$-tight if

\vspace{0.03in} $(i)$ For all $T>0$, there exists $n_0$ such that
\beql{sq1}
 J(X_n, T) \equiv \| \Delta X_{n,i} \|_T \equiv \sup_{0 \le t \le T}{\{|X_n (t) - X_n (t-)|\}} \le b_n
\eeq for all $i$ and $n \ge n_0$, where $b_n$ are real numbers
such that $b_n \ra 0$ as $n \ra \infty$; and

\vspace{0.03in} $(ii)$ The sequence of optional quadratic
variation processes $[ X_{n,i}]$ is $C$-tight for each $i$.
 \end{theorem}

In closing we remark that
 these quadratic-variation tightness conditions in Theorems \ref{thAngle} and \ref{thSquare}
 are often easy to verify because these are nondecreasing processes; see Lemma \ref{lemModQCV}.

\section{Proofs of Tightness}\label{secTightProof}

We now begin the proof of the martingale FCLT in Theorem \ref{thMart}. In this section we
do the tightness part; in the next section we do the
characterization part.

\subsection{Quick Proofs of $C$-Tightness}\label{secFast}

We can give very quick proofs of $C$-tightness in both cases,
exploiting Theorems \ref{thAngle} and \ref{thSquare}. In case (i)
we need to add an extra assumption.  We need to assume that the
jumps are uniformly bounded and that the bound is asymptotically
negligible. In particular, we need to go beyond conditions
\eqn{ek1} and assume that, for all $T>0$, there exists $n_0$ such
that \beql{fast1} J(M_{n,i}, T)  \le b_n \qforallq i \qandq n \ge
n_0 ~, \eeq where $J$ is the maximum-jump function in \eqn{mod3}
and $b_n$ are real numbers such that $b_n \ra 0$ as $n \ra
\infty$, as in Theorem \ref{thSquare}.  We remark that this extra
condition is satisfied for many queueing applications, as in Pang et al. \cite{PTW07},
because the jumps are at most of size $1$ before scaling, so that
they become at most of size $1/\sqrt{n}$ after scaling.

\paragraph{Case (i) with bounded jumps.}

Suppose that the jumps of the martingale are uniformly bounded and
that the bound is asymptotically negligible, as in \eqn{fast1}.
We thus can apply Theorem \ref{thSquare}.  By condition \eqn{ek2}
and Lemma \ref{lemModQCV}, we have convergence of the optional
quadratic-variation (square-bracket) processes:  $[M_{n,i},
M_{n,j}] \Rightarrow c_{i,j} e$ in $D$ for each $i$ and $j$.  Thus
these sequences of processes are necessarily $C$-tight.  Hence
Theorem \ref{thSquare} implies that the sequence of martingales
$\{M_{n,i}: n \ge 1\}$ is $C$-tight in $D$ for each $i$.  Lemma
\ref{lemTightProd} then implies that $\{M_{n}: n \ge 1\}$ is
$C$-tight in $D^k$.~~~\bsq

\paragraph{Case (ii).}

We need no extra condition in case (ii).  We can apply Theorem
\ref{thAngle}. Condition \eqn{ek6} and Lemma \ref{lemModQCV} imply
that there is convergence of the predictable quadratic-variation
(angle-bracket) processes: $\langle M_{n,i}, M_{n,j} \rangle
\Rightarrow c_{i,j} e$ in $D$ for each $i$ and $j$.  Thus these
sequences of processes are necessarily $C$-tight, so condition
(ii) of Theorem \ref{thAngle} is satisfied.  Condition (i) is
satisfied too, since $M_{n, i} (0) = 0$ for all $i$ by
assumption.~~~\bsq

\subsection{The Ethier-Kurtz Proof in Case (i).}\label{secTightness1}

We now give a longer proof for case (i) without making the extra
assumption in \eqn{fast1}. We follow the proof on pp 341-343 of
EK \cite{EK86}, elaborating on several points.

First, if the stochastic processes $M_n$ are local martingales
rather than martingales, then we introduce the available stopping
times $\tau_n$ with $\tau_n \uparrow \infty$ w.p.1 so that $\{M_n
(\tau_n \wedge t): t \ge 0\}$ are martingales for each $n$. Hence
we simply assume that the processes $M_n$ are martingales.

We then exploit the assumed limit $\left[M_{n,i}, M_{n,j}\right]
(t) \Rightarrow c_{i,j} t$ in $\RR$ as $n \ra \infty$ in \eqn{ek2}
in order to act as if $\left[M_{n,i}, M_{n,j}\right] (t)$ is
bounded.  In particular, as in (1.23) on p. 341 of EK, we
introduce the stopping times \beql{ek101} \eta_n \equiv \inf{\{ t
\ge 0: \left[M_{n,i}, M_{n,i}\right] (t) > c_{i,i} t + 1 \quad
\mbox{for some} \quad i\}} ~. \eeq As asserted in EK, by
\eqn{ek2}, $\eta_n \Rightarrow \infty$, but we provide additional
detail about the supporting argument: If we only had the
convergence in distribution for each $t$, it would not follow that
$\eta_n \Rightarrow \infty$. For this step, it is critical that we
can strengthen the mode of convergence to convergence in
distribution in $D$, where the topology on the underlying space
corresponds to uniform convergence over bounded intervals.  To do
so, we use the fact that the quadratic variation processes
$\left[M_{n,i}, M_{n,i}\right] $ are monotone and the limit is
continuous. In particular, we apply Lemma \ref{lemModQCV}. But,
indeed, $\eta_n \Rightarrow \infty$ as claimed.

Hence, it suffices to focus on the martingales $\tilde{M}_n \equiv
\{M_n (\eta_n \wedge  t): t \ge 0\}$. We then reduce the $k$
initial dimensions to $1$ by considering the martingales
\beql{ek102} Y_n \equiv \sum_{i=1}^{k} \theta_i \tilde{M}_{n, i}
\eeq for an arbitrary non-null vector $\theta \equiv (\theta_1,
\ldots , \theta_k)$.  The associated optional quadratic variation
process is \beql{ek103} A_{n, \theta} (t)  \equiv [Y_n] = \sum_{i
= 1}^{k} \sum_{j = 1}^{k} \theta_i \theta_j \left[\tilde{M}_{n,i},
\tilde{M}_{n,j}\right] (t), \quad t \ge 0 ~. \eeq From condition
\eqn{ek2} and Lemma \ref{lemModQCV}, $A_{n, \theta} \Rightarrow
c_{\theta} e$ in $D$, where $e$ is the identity function and
\beql{ek104} c_{\theta} \equiv \sum_{i = 1}^{k} \sum_{j = 1}^{k}
\theta_i \theta_j  c_{i,j} ~. \eeq

Some additional commentary may help here.  The topology on the
space $D([0,\infty), \RR^k)$ where the functions take values in
$\RR^k$ is strictly stronger than the topology on the product
space $D([0, \infty), \RR)^k$, but there is no difference on the
subset of continuous functions, and so this issue really plays no
role here.  We can apply Lemma \ref{lemContLim} and condition
\eqn{ek1} to conclude that a limit of any convergent subsequence
must have continuous paths. In their use of $Y_n$ defined in
\eqn{ek102}, EK work to establish the stronger tightness in
$D([0,\infty), \RR^k)$; sufficiency is covered by problem 3.22 on
p. 153 of EK. On the other hand, for the product topology it would
suffice to apply Lemma \ref{lemTightProd}, but the form of $Y_n$
in \eqn{ek102} covers that as well. The form of $Y_n$ in
\eqn{ek102} is also convenient for characterizing the limiting
process, as we will see.

The idea now is to establish the tightness of the sequence of
martingales $\{Y_n: n \ge 1\}$ in \eqn{ek102}. The approach in EK
is to introduce smooth bounded functions $f : \RR \ra \RR$ and
consider the associated sequence $\{f(Y_n (t)): n \ge 1\}$.  This
step underlies the entire EK book, and so is to be expected. To
quickly see the fundamental role of this step in EK, see Chapters
1 and 4 in EK, e.g., the definition of a full generator in (1.5.5)
on p 24 and the associated martingale property in Proposition
4.1.7 on p 162. This step is closely related to Ito's formula, as
we explain in \S \ref{secAltChar} below. It is naturally
associated with the Markov-process-centric approach in EK \cite{EK86} and
Stroock and Varadhan \cite{SV79}, as opposed to the martingale-centric
approach in Aldous \cite{A78, A89}, Rebolledo \cite{R80} and JS \cite{JS87}.

Consistent with that general strategy, Ethier and Kurtz exploit
tightness criteria involving such smooth bounded functions; e.g.,
see Theorem \ref{thSub2} and Lemma \ref{lemSimple} here. Having
introduced those smooth bounded functions, we establish tightness
 by applying Lemma \ref{lemSimple}
(ii.a). We will exploit the bounds provided by \eqn{ek101}.  We
let the filtrations \textbf{F}$_n \equiv \{\sF_{n,t}: t \ge 0\}$
be the internal filtrations of $Y_n$.    First, the
submartingale-maximum inequality in Lemma \ref{lemSBmax} will be
used to establish the required stochastic boundedness; see (1.34)
on p. 343 of EK. It will then remain to establish the moment
bounds in \eqn{mod8}.

To carry out both these steps, we fix a function $f$ in
$\sC_{c}^{\infty}$, the set of infinitely differentiable functions
$f: \RR \ra \RR$ with compact support (so that the function $f$
and all derivatives are bounded). The strategy is to apply
Taylor's theorem to establish the desired bounds.  With that in
mind, as in (1.27) on p. 341 of EK, we write
\begin{eqnarray}\label{ek105}
&& E\left[f(Y_n (t+s)) - f(Y_n (t))\left| \right. \sF_{n,t} \right] \nonumber \\
&& \quad = E\left[ \sum_{i=0}^{m-1} \left(f(Y_n (t_{i+1})) - f(Y_n
(t_{i})) - f' (Y_n (t_{i})) \xi_{n,i}\right) \left| \right.
\sF_{n,t}\right] ~,
\end{eqnarray}
where $0 \le t = t_0 < t_1 < \cdots < t_m = t + s$ and $\xi_{n,i}
= Y_n (t_{i+1}) - Y_n (t_{i})$. Formula \eqn{ek105} is justified
because $E[\xi_{n,i}| \sF_{n,t}] = 0$ for all $i$ since $Y_n$ is
an $\bf{F_n}$-martingale. Also, there is cancellation in the first
two terms in the sum.

Following (1.28)--(1.30) in EK, we write \beql{ek106} \gamma_n
\equiv \max{\{j: t_j < \eta_n \wedge (t+s)\}} \eeq and
\begin{eqnarray}\label{ek107}
&&\zeta_n \equiv  \\
&& \quad \max{\{j: t_j < \eta_n \wedge (t+s), \sum_{i = 0}^{j}
\xi_{n,i}^2 \le  k \sum_{i=1}^{k} c_{i,i} \theta_i^2 (t+s) + 2k
\sum_{i=1}^k \theta_i^2\}} ~. \nonumber
\end{eqnarray}
By the definition of $\eta_n$, we will have $\gamma_n = \zeta_n$
when there are sufficiently many points so that the largest
interval $t_{i+1} - t_i$ is suitably small and $n$ is sufficiently
large. Some additional explanation might help here: To get the
upper bound in \eqn{ek107}, we exploit
the definition of the optional quadratic covariation (equation (24) in Pang et al. \cite{PTW07})
 and the simple
relation $(a + b)^2 \ge 0$, implying that $|2ab| \le  a^2 + b^2$
for real numbers $a$ and $b$. Thus \beql{ek107a} 2 | [M_{n,i},
M_{n,j}] (t) \theta_i \theta _j | \le  [M_{n,i}, M_{n,i}] (t)
\theta_i^2 +  [M_{n,j}, M_{n,j}] (t) \theta _j^2 \eeq for each $i$
and $j$, so that \beql{ek107b} |A_{n,\theta} (t) |  =  | \sum_{i =
1}^{k} \sum_{j = 1}^{k} \theta_i \theta_j \left[\tilde{M}_{n,i},
\tilde{M}_{n,j}\right] (t) |
 \le   k \sum_{i = 1}^{k}  \left[\tilde{M}_{n,i}, \tilde{M}_{n,i}\right] (t) \theta_i^2 ~.
\eeq

In turn, by \eqn{ek101}, inequality \eqn{ek107b} implies that
\beql{ek107c} |A_{n,\theta} (t) | \le  k \sum_{i = 1}^{k} (c_{i,i}
t + 1) \theta_i^2 =   k \sum_{i = 1}^{k} c_{i,i} t \theta_i^2 +
k \sum_{i = 1}^{k}  \theta_i^2 ~. \eeq In definition \eqn{ek107}
we increase the target by replacing the constant $1$ by $2$ in the
last term, so that eventually we will have $\zeta_n = \gamma_n$,
as claimed.

As in (1.30) on p. 341 of EK, we next extend the expression
\eqn{ek105} to include second derivative terms, simply by adding
and subtracting.  In particular,
 we have
\begin{eqnarray}\label{ek108}
&&\hspace*{-.7cm} E\left[f(Y_n (t+s)) - f(Y_n (t))\left| \right. \sF_{n,t} \right]  \nonumber \\
&& \hspace*{-.7cm} = E\left[ \sum_{i=\zeta_n}^{\gamma_n} \left(f(Y_n (t_{i+1}))
- f(Y_n (t_{i})) - f' (Y_n (t_{i}))\xi_{n,i} \right) \left| \right. \sF_{n,t}\right] \nonumber \\
&&\hspace*{-.7cm} \quad  + E\left[ \sum_{i=0}^{\zeta_n -1} \left(f(Y_n (t_{i+1}))
- f(Y_n (t_{i})) - f' (Y_n (t_{i}))\xi_{n,i} - \frac{1}{2} f'' (Y_n (t_{i}))\xi_{n,i}^2 \right) \left| \right. \sF_{n,t}\right]
\nonumber\!\!\!\! \\
&&\hspace*{-.7cm} \quad + E\left[ \sum_{i=0}^{\zeta_n -1} \frac{1}{2} f'' (Y_n
(t_{i}))\xi_{n,i}^2 \left| \right. \sF_{n,t}\right] ~.
\end{eqnarray}
Following p. 342 of EK, we set $\Delta Y_n (u) \equiv Y_n (u) -
Y_n (u-)$ and introduce more time points so that the largest
difference $t_{i+1} - t_i$ converges to $0$.  In the limit as we
add more time points,
 we obtain the representation
\begin{eqnarray}\label{ek108a}
E\left[f(Y_n (t+s)) - f(Y_n (t))\left| \right. \sF_{n,t} \right] & = & W_{n,1} (t, t+s) + W_{n,2} (t, t+s) \nonumber \\
&& \quad + W_{n,3} (t, t+s) ~,
\end{eqnarray}
where
\begin{eqnarray}\label{ek108b}
\hspace*{-.34cm}W_{n,1} (t, t+s) \hspace*{-.2cm} &\equiv&\hspace*{-.2cm} E \left[f(Y_n ((t+s)\wedge \eta_n))
  - f(Y_n ((t+s)\wedge \eta_n)-) \right. \nonumber \\
\hspace*{-.34cm} && \quad  -  \left. f' (Y_n ((t+s)\wedge \eta_n )-) \Delta Y_n ((t+s)\wedge \eta_n)) \left| \right. \sF_{n,t}  \right] \nonumber \\
\hspace*{-.34cm}W_{n,2} (t, t+s) \hspace*{-.2cm} &\equiv& \hspace*{-.2cm}    E\Biggl[ \sum_{\,t \wedge \eta_n < u <
(t+s)\wedge \eta_n }\!\!\!\! \left( f(Y_n (u))
\,{-}\, f(Y_n (u-)) \,{-}\, f' (Y_n (u-))\Delta Y_n (u) \right. \nonumber\hspace*{-1cm} \\
\hspace*{-.34cm}&& \quad  - f' (Y_n (u-))\Delta Y_n (u)
- \left. \frac{1}{2} f'' (Y_n (u-))(\Delta Y_n (u))^2  \right) \left| \right. \sF_{n,t} \Biggr] \nonumber \\
\hspace*{-.34cm}W_{n,3} (t, t+s) \hspace*{-.2cm}&\equiv&\hspace*{-.2cm}  E\left[ \int_{t \wedge
\eta_n}^{(t+s)\wedge \eta_n} \frac{1}{2} f'' (Y_n (u-))\,
dA_{n,\theta} (u) \left| \right. \sF_{n,t} \right] ~.
\end{eqnarray}
(We remark that we have inserted $t \wedge \eta_n$ in two places
in (1.31) on p 342 of EK. We can write the sum for $W_{n,2}$ that
way because $Y_n$ has at most countably many discontinuities.)

Now we can bound the three terms in \eqn{ek108a} and \eqn{ek108b}.
As in (1.32) on p. 342 of EK, we can apply Taylor's theorem in the
form \beql{Taylor} f(b) = f(a) + f^{(1)}(a) (b-a) + \cdots +
f^{(k-1)} (a) \frac{(b-a)^{k-1}}{(k-1)!} + f^{(k)} (c)
\frac{(b-a)^k}{k!} \eeq for some point $c$ with $a < c < b$ for $k
\ge 1$ (using modified derivative notation).

Applying this to the second conditional-expectation term in
\eqn{ek108a} and the definition of $\eta_n$ in \eqn{ek101}, we get

\beq W_{n,2} (t, t+s) \le   \frac{\|f'''\|}{6}
  E\left[ \sup_{t < u < t+s} | \Delta Y_n (u)|   k \sum_{i=1}^{k} ( c_{i,i} (t+s) + 1) \theta_i^2   \left| \right. \sF_{n,t}\right] ~.
  \eeqno

  Reasoning the same way for the other two terms, overall we get a minor variation of the bound in (1.33) on p. 342 of EK, namely,
\begin{eqnarray}\label{ek110}
&& E\left[f(Y_n (t+s)) - f(Y_n (t))\left| \right. \sF_{n,t} \right]  \nonumber \\
&& \quad \le C_{f} E\left[ \sup_{t < u \le t+s} | \Delta Y_n (u)|
\left( 1 +   k \sum_{i=1}^{k} ( c_{i,i} (t+s) + 1) \theta_i^2  \right) \right. \nonumber \\
&& \quad \quad \left. + A_{n, \theta} ((t+s) \wedge \eta_n -) -
A_{n, \theta} (t \wedge \eta_n -) \left| \right. \sF_{n,t}\right]
~,
\end{eqnarray}
where the constant $C_f$ depends only on the norms of the first
three derivatives:  $\|f'\|$, $\|f''\|$ and $\|f'''\|$.

We now apply the inequality in \eqn{ek110} to construct the
required bounds.  First, for the stochastic boundedness needed in
condition (i) of Lemma \ref{lemSimple}, we let the function $f$
have the additional properties in Lemma \ref{lemSBmax}; i.e., we
assume that $f$ is an even nonnegative convex function with $f'
(t) > 0$ for $t > 0$.  Then there exists a constant $K$ such that
\beql{ek111} E[f(Y_n (T))]  \le  C_{f} E\left[ \left(1 + J(Y_n,
T)|\right)
 \left( 1 +   k \sum_{i=1}^{k} ( c_{i,i} (T) + 1) \theta_i^2  \right) \right] < K
\eeq for all $n$ by virtue of \eqn{ek1}, again using $J$ in
\eqn{mod3}).  (See (1.34) on p 343 of EK, where $t$ on the right
should be $T$.)

Next, we can let the random variables $Z_n (\delta, f, T)$ needed
in condition (ii.a) of Lemma \ref{lemSimple} be defined, as in
(1.36) on p. 343 of EK, by
\begin{eqnarray}\label{ek112}
Z_n (\delta, f, T) & \equiv &   C_f \left[ J(Y_n, T+ \delta)
 \left( 1 +   k \sum_{i=1}^{k} ( c_{i,i} (t+s) + 1) \theta_i^2  \right) \right. \nonumber \\
&& \quad  + \left. \sup_{0 \le t \le T}{\{ A_{n, \theta}
((t+\delta) \wedge \eta_n -) - A_{n, \theta} (t \wedge \eta_n
-)\}}\right] ~,
\end{eqnarray}
where $C_f$ is a new constant depending on the function $f$. Then,
by \eqn{ek1} and \eqn{ek2}, \beql{ek113} \lim_{\delta \downarrow
0} \limsup_{n \ra \infty} E[Z_n (\delta, f, T)] = 0 ~, \eeq as
required in condition (ii.a) of Lemma \ref{lemSimple}. That
completes the proof of tightness in case (i).~~~\bsq

\subsection{The Ethier-Kurtz Proof in Case (ii).}\label{secTightness2}

We already gave a complete proof of tightness for case (ii) in \S
\ref{secFast}. We now give an alternative proof, following EK. The
proof starts out the same as the proof for case (i). First, if the
stochastic processes
$M_{n,i,i}^2 - \langle M_{n,i}, M_{n,i} \rangle$ are local
martingales rather than martingales, then we introduce the
available stopping times $\tau_n$ with $\tau_n \uparrow \infty$
w.p.1 so that
$\{M_{n,i,i}^2 (\tau_n \wedge t) - \langle M_{n,i}, M_{n,i}
\rangle (\tau_n \wedge t): t \ge 0\}$
 are martingales for each $n$ and $i$.  Hence, we can assume that
$M_{n,i,i}^2 - \langle M_{n,i}, M_{n,i} \rangle$ are martingales.

We then exploit the assumed limit $\langle M_{n,i}, M_{n,j}
\rangle (t) \Rightarrow c_{i,j} t$ in $\RR$ as $n \ra \infty$ in
\eqn{ek6} in order to act as if $\langle M_{n,i}, M_{n,j} \rangle
(t)$ is bounded.  In particular, as in \eqn{ek101} and (1.23) on
p. 341 of EK, we introduce the additional stopping times
\beql{ek201} \eta_n \equiv \inf{\{ t \ge 0: \langle M_{n,i},
M_{n,i} \rangle (t) > c_{i,i} t + 1 \quad \mbox{for some} \quad
i\}} ~. \eeq We then define \beql{ek202} \tilde{M}_n (t) \equiv
M_n (\eta_n \wedge t), \quad t \ge 0 ~. \eeq Again we apply Lemma
\ref{lemModQCV} to deduce that we have convergence $\langle
M_{n,i}, M_{n,i} \rangle \Rightarrow c_{i,i} e$ in $D$ for each
$i$, which implies that $\eta_n \Rightarrow \infty$ as $n \ra
\infty$.

Closely paralleling EK, simplify notation by writing \beql{ek203}
A_{n, i,j} (t) \equiv \langle M_{n,i}, M_{n,j} \rangle (t) \qandq
\tilde{A}_{n, i,j} (t) \equiv A_{n, i,j} (t \wedge \eta_n), \quad
t \ge 0 ~. \eeq Then \beql{ek204} \tilde{A}_{n, i,i} (t) \le 1 +
c_{i,i} t + J(A_{n, i,i},t), \quad t \ge 0 ~, \eeq for $J$ in
\eqn{mod3}, where \beql{ek205} E[J(A_{n, i,i},t)] \ra 0 \qasq n
\ra \infty \qforallq t > 0 ~, \eeq by condition \eqn{ek4}.  Since
$\tilde{M}_{n,i,i}^2 - \tilde{A}_{n, i,i} $ is a martingale, so is
$\sum_{i=1}^k (\tilde{M}_{n,i,i}^2 - \tilde{A}_{n, i,i} \rangle)$.
Hence, for $0 \le s < t$, \beql{ek206}
 \sum_{i=1}^k E\left[ \tilde{M}_{n, i} (t +s)^2 - \tilde{M}_{n, i} (t)^2 | \sF_{n,t}\right]
 =  \sum_{i=1}^k  E\left[ \tilde{A}_{n, i} (t +s) - \tilde{A}_{n, i} (t) | \sF_{n,t}\right]
\eeq On the other hand,
\begin{eqnarray}\label{ek207}
&& \sum_{i=1}^k E\left[  \left(\tilde{M}_{n, i} (t +s) - \tilde{M}_{n, i} (t)\right)^2 | \sF_{n,t}\right] \nonumber \\
&& \quad =   \sum_{i=1}^k E\left[\left(\tilde{M}_{n, i} (t +s)^2 -
2 \tilde{M}_{n, i} (t +s) \tilde{M}_{n, i} (t)
+\tilde{M}_{n, i} (t)^2 \right)| \sF_{n,t}\right] \nonumber \\
&& \quad =   \sum_{i=1}^k E \left[ \left( \tilde{M}_{n, i} (t +
s)^2 - \tilde{M}_{n, i} (t)^2 \right) | \sF_{n,t}\right] ~.
\end{eqnarray}
Hence,
\begin{eqnarray}\label{ek208}
&&  E\left[ \sum_{i=1}^k \left(\tilde{M}_{n, i} (t +s) - \tilde{M}_{n, i} (t)\right)^2 | \sF_{n,t}\right]\nonumber \\
 && \quad = E\left[ \sum_{i=1}^k  \tilde{A}_{n, i} (t +s) - \tilde{A}_{n, i} (t) | \sF_{n,t}\right] ~,
\end{eqnarray}
as in 1.41 on p 344 of EK.

Now we can verify condition (ii.b) in Lemma \ref{lemSimple}: We
can define the random variables $Z_n (\delta, T)$ by \beql{ek209}
Z_n (\delta, T) \equiv \sup_{0 \le t \le T}{\{ \sum_{i=1}^k (
\tilde{A}_{n,i, i} (t + \delta) - \tilde{A}_{n,i, i} (t))\}} ~.
\eeq By condition \eqn{ek6} and Lemma \ref{lemModQCV},
\beql{ek209a} Z_n (\delta, T) \Rightarrow \sup_{0 \le t \le T}{\{
\sum_{i=1}^k \left( c_{i,i} (t+ \delta) - c_{i,i} (t)\right)\}} =
\delta \sum_{i=1}^k c_{i,i} \qasq n \ra \infty ~. \eeq We then
have uniform integrability because, by \eqn{ek204} and
\eqn{ek209}, \beql{ek210} Z_n (\delta, T) \le  \sum_{i=1}^k \left(
1 + c_{i,i} (T + \delta) + J(A_{n,i, i}, T)\right) ~, \eeq where
$J$ is the maximum jump function in \eqn{mod3} and $E
\left[J(A_{n,i, i}, T)\right] \ra 0$ as $n \ra \infty$. Hence,
\beql{ek211} \lim_{\delta \downarrow 0} \limsup_{n \ra \infty}
E[Z_n (\delta, T)] =  0 ~. \eeq

Finally, we must verify condition (i) Lemma \ref{lemSimple}; i.e.,
we must show that $\{M_{n,i}: n \ge 1\}$ is stochastically
bounded. As for case (i), we can apply Lemma \ref{lemSBmax} for
this purpose, here using the function $f (x) \equiv x^2$.
Modifying \eqn{ek208}, we have, for all $T>0$, the existence of a
constant $K$ such that \beql{ek212}
 E\left[ \sum_{i=1}^k \tilde{M}_{n, i} (T)^2 \right]
 = E\left[ \sum_{i=1}^k  \tilde{A}_{n, i} (T) \right] \le K
\eeq for all sufficiently large $n$.  That completes the
proof.~~~\bsq

\section{Characterization of the Limit}\label{secCharProof}

\subsection{With Uniformly Bounded Jumps.}\label{secChar}

Since the sequences of martingales $\{Y_n: n \ge 1\}$ and
$\{\tilde{M}_{n}: n \ge 1\}$ are tight, they are relatively
compact by Prohorov's theorem (Theorem \ref{thmProhorov}). We
consider a converging subsequence $\tilde{M}_{n_k} \Rightarrow L$.
It remains to show that $L$ necessarily is
 the claimed $k$-dimensional $(0,C)$-Brownian motion $M$.  By Corollary \ref{corProhorov},
 that will imply that $\tilde{M}_{n} \Rightarrow M$.
   Since $\eta_n \Rightarrow \infty$, it then will follow that $M_n \Rightarrow M$ as well.
   But in this section we will be working with the martingales $\tilde{M}_{n_k}$, for which convergence
   has been established by the previous tightness proofs.

   For this characterization step, we first present a proof based
    upon JS, which requires that the jumps be uniformly bounded
    as an extra condition.  We consider the proof in EK afterwards in the next section, \S \ref{secAltChar}.

   It seems useful to state the desired characterization
   result as a theorem.  Conditions (i) and (ii) below cover the two cases of Theorem \ref{thMart}
   under the extra condition \eqn{extra}.
\begin{theorem}{\em $($characterization of the limit$)$}\label{thChar}
Suppose that $M_n \Rightarrow M$ in $D^k$, where $M_n \equiv
(M_{n,1}, \ldots , M_{n,k})$ is a $k$-dimensional local martingale
adapted to the filtration $\textbf{F}_n \equiv \{\sF_{n,t}\}$ for
each $n$ with bounded jumps, i.e., for all $T > 0$, there exists
$n_0$ and $K$ such that \beql{extra} J(M_{n,i}, T) \le K \qforallq
n \ge n_0 ~, \eeq where $J$ is the maximum-jump function in {\em
\eqn{mod3}}. Suppose that $M \equiv (M_1, \ldots , M_k)$ is a
continuous $k$-dimensional process with $M(0) = (0, \ldots, 0)$.
In addition, suppose that either \beql{ch1}
 (i) \quad [M_{n,i}, M_{n,j}] (t) \Rightarrow c_{i,j} t \qasq n \ra \infty \qinq \RR
\eeq for all $t > 0$, $i$ and $j$, or (ii)   $M_n$ is locally
square integrable and \beql{ch2} \langle M_{n,i}, M_{n,j} \rangle
(t) \Rightarrow c_{i,j} t \qasq n \ra \infty \qinq \RR \eeq for
all $t > 0$, $i$ and $j$. Then $M$ is a $k$-dimensional
$(0,C)$-Brownian motion, having time-dependent mean vector and
covariance matrix \beql{ek8z} E[M(t)] = (0, \ldots , 0) \qandq
E[M(t)M(t)^{tr}] = C t, \quad t \ge 0 ~. \eeq
\end{theorem}

Theorem \ref{thChar} follows directly from several other basic
results.
 The first is the classical
L\'{e}vy martingale characterization of Brownian motion; e.g., see
p. 156 of Karatzas and Shreve \cite{KS88} or p. 102 of JS. Ito's
formula provides an elegant proof.  (Theorem 7.1.2 on p 339 of EK
is a form that exploits Ito's formula in the statement.) Recall
that for a continuous local martingale both the angle-bracket and
square-bracket processes are well defined and equal.

\begin{theorem}{\em $($L\'{e}vy martingale characterization of Brownian motion$)$}\label{thLevy}
Let $M \equiv (M_1, \ldots , M_k)$ be a continuous $k$-dimensional
process adapted to a filtration $\textbf{F} \equiv \{\sF_t\}$ with
$M_i (0) = 0$ for each $i$. Suppose that each one-dimensional
marginal process $M_i$ is a continuous
local-$\textbf{F}$-martingale.  If either the optional covariation
processes satisfy \beql{Levy1} [M_i, M_j] = c_{i,j} t , \quad t
\ge 0 ~, \eeq for each $i$ and $j$, or if $M_i$ are locally
square-integrable martingales, so that the predictable quadratic
covariation processes are well defined, and they satisfy
\beql{Levy2} \langle M_i, M_j \rangle (t) = c_{i,j} t, \quad t \ge
0 ~, \eeq for each $i$ and $j$, where $C \equiv (c_{i,j})$ is a
nonnegative-definite symmetric real matrix, then $M$ is a
$k$-dimensional $(0,C)$-Brownian motion, i.e., \beql{ek8Levy}
E[M(t)] = (0, \ldots , 0) \qandq E[M(t)M(t)^{tr}] = C t, \quad t
\ge 0 ~. \eeq
\end{theorem}

We now present two basic preservation results. The first
preservation result states that the limit of any convergent
sequence of martingales must itself be a martingale, under
regularity conditions.  The following is Problem 7 on p 362 of EK.
The essential idea is conveyed by the bounded case, which is
treated in detail in Proposition IX.1.1 on p 481 of JS.  JS then
go on to show that the boundedness can be replaced by uniform
integrability (UI) and then a bounded-jump condition.  (Note that
this is UI of the sequence $\{M_n (t): n \ge 1\}$ over $n$ for
fixed $t$, as opposed to the customary UI over $t$ for fixed $n$.)

\begin{theorem}{\em $($preservation of the martingale property under UI$)$}\label{thPresMG1}
Suppose that $(i)$ $X_n$ is a random element of $D^k$ and $M_n$ is
a random element of $D$ for each $n \ge 1$, $(ii)$ $M_n$ is a
martingale with respect to the filtration generated by $(X_n,
M_n)$ for each $n \ge 1$, and $(iii)$ $(X_n, M_n) \Rightarrow
(X,M)$ in $D^{k+1}$ as $n \ra \infty$. If, in addition, $\{M_n
(t): t \ge 1\}$ is uniformly integrable for each $t > 0$, then $M$
is a martingale with respect to the filtration generated by
$(X,M)$ $($and thus also with respect to the filtration generated
by $M)$.
\end{theorem}

\paragraph{Proof.}  Let $\{ \sF_{n,t}: t \ge 0\}$ and $\{ \sF_{t}: t \ge 0\}$ denote the filtrations generated by $(X_n, M_n)$
and $(X, M)$, respectively, on their underlying probability
spaces.  Let the probability space for the limit $(X,M)$ be
$(\Omega, \sF, P)$.  It is difficult to relate these filtrations
directly, so we do so indirectly. This involves a rather tricky
measurability argument. We accomplish this goal by considering the
stochastic processes $(X,M)$ and $(X_n,M_n)$ as maps from the
underlying probability spaces to the function space $D^{k+1}$ with
its usual sigma-field generated by the coordinate projections,
i.e., with the filtration \textbf{D$^{k+1}$} $\equiv
\{\sD^{k+1}_t: t \ge0\}$. As discussed in \S VI.1.1 of JS and \S
11.5.3 of Whitt \cite{W02}, the Borel $\sigma$-field on $D^{k+1}$
coincides with the $\sigma$-field on $D^{k+1}$ generated by the
coordinate projections. A $\sD^{k+1}_t$-measurable real-valued
function $f(x)$ defined on $D^{k+1}$ depends on the function $x$
in $D^{k+1}$ only through its behavior in the initial time
interval $[0,t]$.

As in step (a) of the proof of Proposition IX.1.1 of JS, we
consider the map $(X,M): \Omega \ra D^{k+1}$ mapping the
underlying probability space into $D^{k+1}$ with the sigma-field
\textbf{D$^{k+1}$} generated by the coordinate projections.  By
this step, we are effectively lifting the underlying probability
space to $D^{k+1}$, with $(X, M) (t)$ obtained as the coordinate
projection.

Let $t_1$ and $t_2$ with $t_1 < t_2$ be almost-sure continuity
points of the limit process $(X, M)$. Let $f:D^{k+1} \ra \RR$ be a
bounded continuous $\sD^{k+1}_{t_1}$-measurable real-valued
function in that setting.  By this construction, not only is $f(X,
M)$ $\sF_{t_1}$-measurable but $f(X_n, M_n)$ is $\sF_{n, t_1}$
measurable for each $n \ge 1$ as well.  (Note that $f(X_n, M_n)$
is the composition of the maps $(X_n, M_n) : \Omega_n \ra D^{k+1}$
and $f: D^{k+1} \ra \RR$.)

By the continuous-mapping theorem, we have first \beq \left(f(X_n,
M_n), M_n (t_2), M_n (t_1)\right) \Rightarrow \left(f(X,M), M
(t_2), M (t_1)\right) \qinq  \RR^3 \qasq n \ra \infty \eeqno and
then \beq f(X_n, M_n)(M_n (t_2)- M_n (t_1)) \Rightarrow f(X,M) ( M
(t_2) - M (t_1)) \qinq  \RR \qasq n \ra \infty ~. \eeqno By the
boundedness of $f$ and the uniform integrability of $\{M_n
(t_i)\}$, we then have \beq E\left[f(X_n, M_n)(M_n (t_2)- M_n
(t_1))\right] \ra E\left[f(X,M) ( M (t_2)- M (t_1))\right]  \qasq
n \ra \infty ~. \eeqno Now we exploit the fact that $f(X_n, M_n)$
is actually $\sF_{n, t_1}$-measurable for each $n$. We are thus
able to invoke the martingale property for each $n$ and conclude
that \beq E\left[f(X_n, M_n)(M_n (t_2)- M_n (t_1))\right] = 0
\qforallq n ~. \eeqno Combining these last two relations, we have
the relation \beql{pr201} E\left[f(X,M) ( M (t_2)- M (t_1))\right]
= 0 ~. \eeq Now, by a monotone class argument (e.g., see p 496 of
EK), we can show that the relation \eqn{pr201} remains true for
all bounded $\sD^{k+1}_{t_1}$-measurable real-valued functions
$f$, which includes the indicator function of an arbitrary
measurable set $A$ in $\sD^{k+1}_{t_1}$. Hence, \beq
E\left[1_{\{(X,M)\in A\}} ( M (t_2)- M (t_1))\right] = 0 ~. \eeqno
This in turn implies that \beq E\left[1_{B} ( M (t_2)- M
(t_1))\right] = 0 \eeqno for all $B$ in $\sF_{t_1}$, which implies
that \beq E\left[M (t_2) - M (t_1))| \sF_{t_1} \right] = 0 ~,
\eeqno which is the desired martingale property for these special
time points $t_1$ and $t_2$. We obtain the result for arbitrary
time points $t_1$ and $t_2$ by considering limits from the right.
~~~\bsq

A way to get the uniform-integrability regularity condition for
martingales in Theorem \ref{thPresMG1} is
to have a local martingale with uniformly bounded jumps. The
following is Corollary IX.1.19 to Proposition IX.1.17 on p 485 of
JS.

\begin{theorem}{\em $($preservation of the local-martingale property with bounded\break jumps$)$}\label{thPresMG2}
Suppose that conditions $(i)-(iii)$ of Theorem {\em
\ref{thPresMG1}} hold, except that $M_n$ is only required to be a
local martingale for each $n$. If, in addition,
 for each $T > 0$, there exists
a positive integer $n_0$ and a constant $K$ such that \beql{pres1}
J(M_n, T)  \le K \qforallq n \ge n_0 ~, \eeq then $M$ is a local
martingale with respect to the filtration generated by $X$ $($and
thus also with respect to the filtration generated by $M)$.
\end{theorem}

The second preservation result states that, under regularity
conditions, the optional quadratic variation $[M]$ of a local
martingale is a continuous function of the local martingale.  The
following is Corollary VI.6.7 on p 342 of JS.

\begin{theorem}{\em $($preservation of the optional quadratic variation$)$}\label{thPresOQV}
Suppose that $M_n$ is a local martingale for each $n \ge 1$ and
$M_n \Rightarrow M$ in $D$. If, in addition, for each $T > 0$,
there exists a positive integer $n_0$ and a constant $K$ such that
\beql{pres2} E\left[J(M_n, T) \right] \le K \qforallq n \ge n_0 ~,
\eeq then \beql{pres3} (M_n, [M_n]) \Rightarrow (M, [M]) \qinq D^2
~. \eeq
\end{theorem}

Note that condition \eqn{pres2} is implied by condition
\eqn{pres1}.  More importantly, condition \eqn{pres2} is implied
by condition \eqn{ek1}.

\paragraph{Proof of Theorem \ref{thChar}.}

From Theorem \ref{thLevy}, we see that it suffices to focus on one
coordinate at a time. To characterize the covariation processes,
we can consider the weighted sums \beql{wt1} M_{n,\theta} \equiv
\sum_{i=1}^{k} \theta_i M_{n,i} \eeq for arbitrary non-null vector
$\theta \equiv (\theta_1, \ldots , \theta_k)$. First suppose that
condition (i) of Theorem \ref{thChar} is satisfied.
 We get the covariations via
\beql{wt2} 2 [M_{n,i}, M_{n,j}] = [M_{n,i} + M_{n,j}, M_{n,i} +
M_{n,j}] - [M_{n,i}, M_{n,i}] - [M_{n,j}, M_{n,j}] ~. \eeq
Henceforth fix $\theta$. First, condition (i) of Theorem
\ref{thChar} implies that \beql{wt3} [ M_{n,\theta}] (t)
\Rightarrow c_{\theta} t \qasq n \ra \infty \qinq \RR \qforallq t
> 0 ~. \eeq Condition \eqn{extra} implies
 condition \eqn{pres1} for $M_{n,\theta}$ .
Hence, we can apply Theorems \ref{thPresMG2} and \ref{thPresOQV}
to deduce that the limit process $M_{\theta} \equiv \sum_{i=1}^{k}
\theta_i M_{i}$ is a local martingale with optional quadratic
variation
 $[M_{\theta}] = c_{\theta} e$, where $c_{\theta} \equiv \sum_{i=1}^{k} \sum_{j= 1}^{k} \theta_i \theta_j c_{i,j}$.
Since this is true for all $\theta$, we can apply Theorem
\ref{thLevy} and \eqn{wt2} to complete the proof.

Now, instead, assume that condition (ii) of Theorem \ref{thChar}
is satisfied. First, by the Doob-Meyer decomposition, Theorem 3.1 of Pang et al. \cite{PTW07}, we can deduce
that $M_{n, \theta}^2 - \langle M_{n, \theta}\rangle$ is a local
martingale for each $n$ and $\theta$. Next,
 it follows from the assumed convergence in \eqn{ch2}, Lemma \ref{lemModQCV} and the continuous mapping theorem
that \beql{wt4} M_{n, \theta}^2 - \langle M_{n, \theta}\rangle
\Rightarrow M_{\theta}^2 - c_{\theta} e \qinq D \qasq n \ra \infty
~. \eeq By Theorem \ref{thPresMG2} and condition \eqn{extra}, the
stochastic process $M_{\theta}^2 - c_{\theta} e$ is a local
martingale for each vector $\theta$, but that means that $\langle
M_{\theta}\rangle (t) = c_{\theta} t$ for each vector $\theta$,
which implies the predictable quadratic covariation condition in
Theorem \ref{thLevy}.  Paralleling \eqn{wt2}, we use \beql{wt5} 2
\langle M_{i}, M_{j} \rangle = \langle M_{i} + M_{j}, M_{i} +
M_{j} \rangle - \langle M_{i}, M_{i} \rangle - \langle M_{j},
M_{j} \rangle ~. \eeq Hence, we can apply Theorem \ref{thLevy} to
complete the proof.~~~\bsq

\subsection{The Ethier-Kurtz Proof in Case (i).}\label{secAltChar}
\hsp The proof of the characterization step for case (i) on p
343-344 of EK is quite brief.  From the pointer to Problem 7 on p.
362 of EK, one might naturally think that we should be applying
the result there, which corresponds to Theorem \ref{thPresMG1}
here, and represent the limit process as a limit of an appropriate
sequence of martingales. A helpful initial observation here is
that it is not actually necessary for the converging processes to
be martingales. Instead of applying Theorem \ref{thPresMG1}, we
can follow the proof of Theorem \ref{thPresMG1} in order to obtain
the desired martingale property asymptotically.

Accordingly, we directly establish that the limit process is a
martingale. In particular, we show that the stochastic process
\beql{xx1} \{f(L(t)) - \frac{c_{\theta}}{2} \int_{0}^{t} f''
(L(s-)) \, ds: t \ge 0\} ~, \eeq is a martingale, where $L$ is the
limit of the converging subsequence $\{Y_{n_k}\}$.
 In particular, we will show that
\beql{xx2} E\left[f(L(t+s)) - f(L(t) - \frac{c_{\theta}}{2}
\int_{t}^{t+s} f'' (L(u-)) \, d u \left| \right. \sF_t \right] = 0
\eeq for $0 < t < t + s$.

Toward that end, EK present their (1.38) and (1.39) on p 343. The
claim about convergence in probability uniformly over compact
subsets of $D$ in (1.38) on p 343 of EK seems hard to verify, but
it is not hard to prove the desired (1.39), exploiting stronger
properties that hold in this particular situation. We now provide
additional details about the proof of a variant of (1.39) on p.
343 of EK. The desired variant of that statement is \eqn{QL2}
below.

\begin{lemma}\label{lemClaim139}
Under the conditions of Theorem {\em \ref{thMart}} in case $(i)$,
if $L$ is the limit of a convergent subsequence $\{Y_{n_k}: n \ge
1\}$, then \beql{QL1} \int_{t\wedge \eta_{n_k}}^{(t+s)\wedge
\eta_{n_k}} \frac{1}{2} f'' (Y_{n_k} (u-)) \, d A_{n_k,\theta} (u)
\Rightarrow \int_{t}^{t+s} \frac{1}{2} f'' (L(u-)) \, c_{\theta}
\, du \eeq in $D([t,t+s], \RR)$ as $n_k \ra \infty$.  As a
consequence, \beql{QL2} E \left[ \left| \right. \int_{t \wedge
\eta_{n_k}}^{(t+s) \wedge \eta_{n_k}} \frac{1}{2} f'' (Y_{n_k}
(u-)) \, d A_{n_k,\theta} (u) - \int_{t}^{t+s} \frac{1}{2} f''
(L(u-)) \, c_{\theta} \, du \left| \right. \right]  \ra 0 ~, \eeq
as in {\em (1.39)} on p {\em 343} of {\em EK}.
\end{lemma}

\paragraph{Proof of Lemma \ref{lemClaim139}.}   We start by proving \eqn{QL1}.
First, we use the fact that both $A_{n,\theta}$ and $c_{\theta} e$
are nondecreasing functions. We also exploit Lemma \ref{lemModQCV}
in order to get $A_{n_k,\theta} \Rightarrow c_{\theta} e$ in $D$.
The topology is uniform convergence over bounded intervals,
because $c_{\theta} e$ is a continuous function.  Given that
$Y_{n_k} \Rightarrow L$, we can apply Lemma \ref{lemContLim} and
condition \eqn{ek1} to deduce that $L$ and thus $f'' (L(t))$
 actually have continuous paths, so the topology is again uniform convergence over bounded intervals.
 Then the limit in \eqn{QL1} is elementary:  First, for any $\epsilon > 0$, $P(\|Y_{n_k} - L\|_{t+s} > \epsilon) \ra 0$
 and $P(\eta_{n_k} \le t+s) \ra 0$
 as $n_k \ra \infty$.  As a consequence, $P(\|f'' \circ Y_{n_k} - f'' \circ L\|_{t+s} > \epsilon) \ra 0$
 as $n \ra \infty$ too.
 Given that $\|f'' \circ Y_{n_k} - f'' \circ L\|_{t+s} \le \epsilon$ and $\eta_{n_k} \ge t + s$,
\begin{eqnarray}\label{QL2a}
 && \left| \int_{t\wedge \eta_{n_k}}^{(t+s)\wedge \eta_{n_k}} \frac{1}{2} f'' (Y_{n_k} (u-)) \, d A_{{n_k},\theta} (u) -
\int_{t}^{t+s} \frac{1}{2} f'' (L(u-)) \, c_{\theta} \, du \right| \nonumber \\
& & \quad \le \left| \int_{t}^{(t+s)} \frac{1}{2} f'' (Y_{n_k}
(u-)) \, d A_{{n_k} ,\theta} (u) -
\int_{t}^{t+s} \frac{1}{2} f'' (L(u-)) \,d A_{n_k,\theta} (u) \right| \nonumber \\
&& \quad + \left| \int_{t}^{t+s} \frac{1}{2} f'' (L(u-)) \, d
A_{{n_k},\theta} (u) -
\int_{t}^{t+s} \frac{1}{2} f'' (L(u-)) \, c_{\theta} \, du \right| \nonumber \\
&& \le  \epsilon ( A_{{n_k},\theta} (t+s) -  A_{{n_k},\theta} (t))
 +  \left| \int_{t}^{t+s} | \frac{1}{2} f'' (L(u-)) | \, ( d A_{{n_k},\theta} (u) -
 c_{\theta} \, du) \right| \nonumber \\
&& \quad \Rightarrow 0 \qasq n \ra \infty ~.
\end{eqnarray}
Hence we have the first limit \eqn{QL1}.
 Since $f''$ and the length of the interval $[t, t+s]$ are both bounded, we also get the associated limit of expectations in \eqn{QL2}.~~~\bsq

\paragraph{Back to the Characterization Proof.}

We now apply the argument used in the proof of Theorem
\ref{thPresMG1}. Let $\{ \sF_{n_k,t}: t \ge 0\}$ and $\{ \sF_{t}:
t \ge 0\}$ denote the filtrations generated by $Y_{n_k}$ and $L$,
respectively, on their underlying probability spaces.  Let the
probability space for the limit $Y$ be $(\Omega, \sF, P)$.  As in
the proof of Theorem \ref{thPresMG1}, we consider a real-valued
function from the function space $D$ with its usual sigma-field
generated by the coordinate projections, i.e., with the filtration
\textbf{D} $\equiv \{\sD_t: t \ge0\}$. A $\sD_t$-measurable
function defined on $D$ depends on the function $x$ in $D$ only
through its behavior in the initial time interval $[0,t]$.

Let $t$ and $t+s$ be the arbitrary time points with $0  \le t < t
+ s$. Let $h:D \ra \RR$ be a bounded continuous
$\sD_{t}$-measurable real-valued function in that setting.  By
this construction, not only is $h(L)$ $\sF_{t}$-measurable, but
$h(Y_{n_k})$ is $\sF_{n_k, t}$ measurable for each $k \ge 1$ as
well.

By the continuous-mapping theorem, we have first
\begin{eqnarray}\label{xx3}
&& \left(h(Y_{n_k}), f(Y_{n_k} (t)), f(Y_{n_k} (t+ s)),
\int_{t\wedge \eta_{n_k}}^{(t+s)\wedge \eta_{n_k}} \frac{1}{2} f'' (Y_{n_k} (u-)) \, d A_{n_k,\theta} (u)  \right) \nonumber \\
&& \quad \Rightarrow \left(h(L), f(L(t)), f(L(t+s)),
 \int_{t}^{t+s} \frac{1}{2} f'' (L(u-)) \, c_{\theta} \, du   \right)
\end{eqnarray}
in $\RR^4$ as $n \ra \infty$, and then
\begin{eqnarray}\label{xx4}
&& h(Y_{n_k}) \left( f(Y_{n_k} (t+ s) -  f(Y_{n_k} (t))
-  \int_{t\wedge \eta_{n_k}}^{(t+s)\wedge \eta_{n_k}} \frac{1}{2} f'' (Y_{n_k} (u-)) \, d A_{n_k,\theta} (u)  \right) \nonumber \\
&& \quad    \Rightarrow h(L) \left( f(L(t+s)) - f(L(t)) -
\int_{t}^{t+s} \frac{1}{2} f'' (L(u-)) \, c_{\theta} \, du
\right)
 \end{eqnarray}
$in$ $\RR$ as $k \ra \infty$. On the other hand, as argued by EK,
by condition \eqn{ek107c} and the boundedness of $f$ and $f''$,
the first two terms $W_{n,1} (t+s)$ and $W_{n,2} (t, t+s)$ on the
right in \eqn{ek108a} converge to $0$ in $L_1$ as well. Hence, the
limit in \eqn{xx4} must be $0$. As a consequence, we have shown
that \beql{xx5} E\left[h(L) \left( f(L(t+s)) - f(L(t)) -
\int_{t}^{t+s} \frac{1}{2} f'' (L(u-)) \, c_{\theta} \, du
\right)\right] = 0 \eeq for all continuous bounded
$\sD_{t}$-measurable real-valued functions $h$.  By the
approximation argument involving the monotone class theorem (as
described in the proof of Theorem \ref{thPresMG1}), we thus get
\beql{xx5} E\left[1_{B} \left( f(L(t+s)) - f(L(t)) -
\int_{t}^{t+s} \frac{1}{2} f'' (L(u-)) \, c_{\theta} \, du
\right)\right] = 0 \eeq for each measurable subset $B$ in
$\sF_{t}$.  That implies \eqn{xx2}, which in turn implies the
desired martingale property. Finally, Theorem 7.1.2 on p 339 of EK
(playing the role of Theorem \ref{thLevy} here) shows that this
martingale property for all these smooth functions $f$ implies
that the limit $L$ must be Brownian motion, as claimed.~~~\bsq

\subsection{The Ethier-Kurtz Proof in Case (ii).}\label{secAltChar2}
\hsp The proof of characterization in this second case seems much
easier for two reasons: First, the required argument is relatively
straightforward (compared to the previous section) and, second, EK
provides almost all the details.

We mostly just repeat the argument in EK. It suffices to focus on
a limit $L$ of a convergent subsequence $\{\tilde{M}_{n_k}: n_k
\ge 1\}$ of the relatively compact sequence $\{\tilde{M}_n: n \ge
1\}$ in $D^k$. By condition \eqn{ek5} and Lemma \ref{lemContLim},
$L$ almost surely has continuous paths.  By the arguments in the
EK tightness proof in \S \ref{secTightness2}, starting with
\eqn{ek206}, we have \beq \sup_{n \ge 1} E[\sum_{i=1}^{k}
\tilde{M}_{n, i} (T)^2] < \infty \eeqno for each $T > 0$.  As a
consequence, the sequence of random vectors $\{\tilde{M}_n (T): n
\ge 1\}$ is uniformly integrable; see (3.18) on p 31 of
Billingsley \cite{B99}.  Hence, Theorem \ref{thPresMG1} implies that
the limit $L$ is a martingale. Moreover, by condition \eqn{ek4}
and Lemma \ref{lemModQCV}, we have the convergence \beq
\tilde{M}_{n_k, i} \tilde{M}_{n_k, j} - \tilde{A}_{n,i,j}
\Rightarrow L_i L_j - c_{i,j} e \qinq D \qasq n_k \ra \infty ~.
\eeqno We can conclude that the limit $L_i L_j - c_{i,j} e$ too is
a martingale by applying Theorem \ref{thPresMG1}, provided that we
can show that the sequence $\{ \tilde{M}_{n_k, i} (T)
\tilde{M}_{n_k, j} (T)  - \tilde{A}_{n_k,i,j} (T): n_k \ge 1\}$ is
uniformly integrable for each $T > 0$.

By \eqn{ek204}, condition \eqn{ek4}, and the inequality \beq 2 |
A_{n_k, i, j} (T) | \le | A_{n_k, i, i} (T) |   + | A_{n_k, j, j}
(T) | ~, \eeqno the sequence $\{ \tilde{A}_{n_k,i,j} (T): n_k \ge
1\}$ is uniformly integrable. Hence it suffices to focus on the
other sequence $\{ \tilde{M}_{n_k, i} (T) \tilde{M}_{n_k, j} (T):
n_k \ge 1\}$, and since \beq 2 | \tilde{M}_{n_k, i} (T)
\tilde{M}_{n_k, j} (T)| \le | \tilde{M}_{n_k, i} (T)^2 |   + |
\tilde{M}_{n_k, j} (T)^2 | ~, \eeqno it suffices to consider the
sequence $\{ \tilde{M}_{n_k, i} (T)^2: n_k \ge 1\}$ for arbitrary
$i$. The required argument now is simplified because the random
variables are nonnegative. Since $\tilde{M}_{n_k, i} (T)^2
\Rightarrow L_i (T)^2$ as $n_k \ra \infty$, it suffices to show
that $E[\tilde{M}_{n_k, i} (T)^2] \ra E[L_i (T)^2]$ as $n_k \ra
\infty$. By \eqn{ek206}--\eqn{ek208}, we already know that
$E[\tilde{M}_{n_k, i} (T)^2] \ra c_{i,i} T$ as $n_k \ra \infty$.
So what is left to prove is that \beql{fin0} E[L_i (T)^2] =
c_{i,i} T ~, \eeq as claimed in (1.45) on p 344 of EK.

For this last step, we introduce stopping times \beql{fin1}
\tau_{n} (x) \equiv \inf{\{ t \ge 0: \tilde{M}_{n, i} (t)^2 > x\}}
\qandq \tau (x) \equiv \inf{\{ t \ge 0: L_{i} (t)^2 > x\}} ~, \eeq
for each $x > 0$. We deviate from EK a bit here.  It will be
convenient to guarantee that these stopping times are finite for
all $x$. (So far, we could have $L(t) = 0$ for all $t \ge 0$.)
Accordingly, we define modified stopping times that necessarily
are finite for all $x$. We do so by adding the continuous
deterministic function $(t - T)^{+}$ to the two stochastic
processes. In particular, let \beq \tau_{n} (x)' \equiv \inf{\{ t
\ge 0: \tilde{M}_{n, i} (t)^2 + (t - T)^{+} > x\}} \eeqno and \beq
\tau (x)' \equiv \inf{\{ t \ge 0: L_{i} (t)^2 + (t - T)^{+} > x\}}
\qforq x \ge 0~. \eeqno With this modification, the stopping times
$\tau_{n} (x)'$ and $\tau (x)'$ are necessarily finite for all
$x$, and yet nothing is changed over the relevant time interval
$[0,T]$: \beql{fin2a} \tau_{n} (x)' \wedge t = \tau_{n} (x) \wedge
t \qandq \tau (x)' \wedge t = \tau (x) \wedge t \qforq 0 \le t \le
T \eeq for all $x \ge 0$.

After having made this minor modification, we can apply the
continuous mapping theorem with the inverse function as in \S 13.6
of Whitt \cite{W02}, mapping the subset of functions unbounded above
in $D$ into itself, which for the cognoscenti will be a
simplification.  That is, we regard the first-passage times as
stochastic processes indexed by $x \ge 0$. For this step, we work
with the stochastic processes (random elements of $D$) $\tau_{n}
\equiv \{\tau_{n} (x): x \ge 0\}$ and similarly for the other
first passage times. We need to be careful to make this work:  We
need to use a weaker topology on D on the range. But Theorem
13.6.2 of \cite{W02} implies that \beql{fin3a} (\tau_{n_k}',
\tilde{M}_{n_k}) \Rightarrow (\tau', L) \qinq (D, M_{1}^{'})
\times (D, J_1) \qasq k \ra \infty ~, \eeq which in turn, by the
continuous mapping theorem again, implies that \beql{fin3b}
(\tau_{n_k} (x)' \wedge T, \tilde{M}_{n_k}) \Rightarrow (\tau (x)'
\wedge T, L) \qinq  \RR \times (D, J_1) \qasq k \ra \infty ~, \eeq
for all $x$ except the at most countably many $x$ that are
discontinuity points of the limiting stochastic process $\{\tau
(x)': x \ge 0\}$. Since we have restricted the time argument to
the interval $[0,T]$, we can apply \eqn{fin2a} and use the
original stopping times, obtaining \beql{fin3b} (\tau_{n_k} (x)
\wedge T, \tilde{M}_{n_k}) \Rightarrow (\tau (x) \wedge T, L)
\qinq  \RR \times (D, J_1) \qasq k \ra \infty \eeq again for all
$x$ except the at most countably many $x$ that are discontinuity
points of the limiting stochastic process $\tau' \equiv \{\tau
(x)': x \ge 0\}$. By the continuous mapping theorem once more,
with the composition map and the square, e.g., see VI.2.1 (b5) on
p 301 of JS, we obtain \beql{fin3c} \tilde{M}_{n_k, i} (t \wedge
\tau_{n_k} (x))^2 \Rightarrow  L (t \wedge \tau (x))^2 \qinq \RR
\qasq k \ra \infty \eeq for all but countably many $x$.  At this
point we have arrived at a variant of the conclusion reached on p
345 of EK. (We do not need to restrict the $T$ because the limit
process $L$ has been shown to have continuous paths.) It is not
necessary to consider the inverse function mapping a subset of $D$
into itself, as we did; instead we could go directly to
\eqn{fin3b}, but that joint convergence is needed in order to
apply the continuous mapping theorem to get \eqn{fin3c}.

Now we have the remaining UI argument in (1.48) of EK, which
closely parallels \eqn{ek210}.  We observe that the sequence
$\{\tilde{M}_{n_k,i} (T \wedge \tau_{n_k} (x))^2: k \ge 1\}$ is
UI, because \beql{fin3d} \tilde{M}_{n_k,i} (T \wedge
\tau_{n_k}(x))^2 \le 2 \left( x + J(\tilde{M}_{n_k,i}^2,T)
\right) \eeq where \beql{fin3e} E\left[  J(\tilde{M}_{n_k,i}^2,T)
\right] \ra 0 \qasq k \ra \infty \eeq by virtue of condition
\eqn{ek5}, again using $J$ in \eqn{mod3}.

As a consequence, as in (1.49) on p 345 of EK, for all but
countably many $x$, we get
\begin{eqnarray}\label{fin5}
E[L_i (T \wedge \tau (x))^2] & = & \lim_{k \ra \infty} E[\tilde{M}_{n_k, i} (T \wedge \tau_{n_k} (x))^2] \nonumber \\
& = & \lim_{k \ra \infty} E[\tilde{A}_{n_k, i, i} (T \wedge \tau_{n_k} (x))] \nonumber \\
& = & E[ c_{i,i} (T \wedge \tau (x))] ~.
\end{eqnarray}
Letting $x \ra \infty$ through allowed values, we get $\tau (x)
\Rightarrow \infty$ and then the desired \eqn{fin0}.~~~\bsq

\section*{Acknowledgments.}  The author thanks Tom Kurtz and Columbia doctoral students
 Itay Gurvich, Guodong Pang and Rishi Talreja for helpful comments on this paper.
 The author was supported by
NSF Grant DMI-0457095.


\end{document}